\documentclass{amsart}
\usepackage{amsfonts,amsthm,amssymb,euscript}
\usepackage{graphics}
\usepackage[matrix,arrow]{xy}

\newtheoremstyle{my}{1.5em}{0.5em}{\em}{}{\sc}{:}{0.5em}{}

\theoremstyle{my}
\newtheorem{thm}{Theorem}
\numberwithin{thm}{section}
\numberwithin{equation}{section}
\newtheorem{theorem}[thm]{Theorem}
\newtheorem*{theorem*}{Theorem}
\newtheorem{cor}[thm]{Corollary}

\newtheorem*{corollary*}{Corollary}
\newtheorem{lemma}[thm]{Lemma}

\newtheorem{prop}[thm]{Proposition}

\newtheorem{conjecture}[thm]{Conjecture}
\newtheorem*{conjecture*}{Conjecture}

\newtheorem*{question*}{Question}

\newtheorem{definition}[thm]{Definition}

\newtheorem*{definitions*}{Definitions}

\newtheorem*{rem*}{Remark}
\newtheorem{remark}[thm]{Remark}

\newtheorem*{remark*}{Remark}

\newtheorem*{remarks*}{Remarks}
\newtheorem*{example*}{Example}
\newtheorem{example}[thm]{Example}
\newtheorem*{examples*}{Examples}

\newcommand{\R}{\mathbb{R}}
\newcommand{\Z}{\mathbb{Z}}
\newcommand{\Q}{\mathbb{Q}}
\newcommand{\C}{\mathbb{C}}

\newcommand{\half}{{\textstyle\frac{1}{2}}}
\newcommand{\quarter}{{\textstyle\frac{1}{4}}}

\newcommand{\iso}{\cong}           
\newcommand{\smooth}{C^\infty}
\newcommand{\CP}[1]{\C {\mathrm P}^{#1}}

\newcommand{\leftsc}{\langle}
\newcommand{\rightsc}{\rangle}

\renewcommand{\o}{\omega}

\renewcommand{\subsection}[1]{\hspace{-\parindent}\refstepcounter{subsection}{\bf
(\arabic{section}\alph{subsection}) #1:}}

\newcommand{\scrL}{\mathcal{L}}
\newcommand{\scrK}{\mathcal{K}}
\newcommand{\scrM}{\mathcal{M}}
\newcommand{\K}{\mathbb{K}}

\newcommand{\F}{\mathcal{F}}
\newcommand{\maurercartan}{\mathcal{MC}}

\title{A biased view of symplectic cohomology}
\author{Paul Seidel}
\date{April 13, 2007}
\begin{document}
\maketitle
\tableofcontents
\newpage

\section{Introduction}

Symplectic cohomology is an invariant of a certain kind of symplectic manifolds (open, or with boundary). It is comparatively easy to define, being a variation on classical Hamiltonian Floer homology. Moreover, its behaviour reflects important aspects of symplectic topology in a fairly direct way. For instance, this applies to the fundamentally trivial nature of subcritical Stein manifolds, and to the importance of exact Lagrangian submanifolds, which are reflected in (nontrivial) vanishing resp.\ non-vanishing theorems for symplectic cohomology. In spite of this, and of the many successful early applications, the theory has not received the same level of attention as, say, Gromov-Witten theory or SFT (symplectic field theory).

These lecture notes are an attempt to advertise the breadth and attractiveness of symplectic cohomology, by stressing connections with various parts of symplectic topology and algebraic geometry. Because of this specific aim, our account may appear somewhat unbalanced (whence the title). On one hand, it includes a certain amount of previously unpublished material (various parts of this are due to Mark McLean, Ivan Smith, and the author; I have tried to mark clearly those places where I am borrowing other people's work). On the other hand, the exposition omits many technical details, and important classical results are stated entirely without proof. Fortunately, there are other surveys which perform much better in these respects, for instance \cite{oancea04b,weber06}. There is also some very recent work which could not be included in these notes, such as that of Bourgeois-Oancea relating symplectic cohomology with contact homology \cite{bourgeois-oancea07a,bourgeois-oancea07b}. Again, others will make up for this deficiency (Cieliebak and Oancea are preparing a paper which will explain this and other SFT-inspired work on symplectic cohomology).

My present understanding of symplectic cohomology has developed throughout the course of many conversations with Ivan Smith. Mohammed Abouzaid, Kevin Costello and Tim Perutz contributed stimulating ideas and comments. I would also like to thank the Harvard mathematics department for inviting me to talk at the 2006 ``Current Developments in Mathematics'' conference, and the audience in my subsequent MIT graduate course for suggestions and corrections. The preparation of these notes was partially funded by NSF grant DMS-0405516.

\section{Liouville domains}

We begin by fixing the class of symplectic manifolds to be used throughout these notes. A precise name might be ``exact symplectic manifolds with contact type boundary'', but for brevity, we will call them Liouville domains. Unlike the case of closed symplectic manifolds, there are several possible notions of isomorphism, with widely different implications. Our definition of ``Liouville isomorphism'' involves attaching an infinite cone, and rules out ``quantitative'' invariants such as volume and capacities.

Generally speaking, constructing Liouville domains is easy because there are no constraints, meaning that the manifold does not have to close up. For instance, if one uses Lefschetz fibrations as a construction method, the global monodromy may be nontrivial, hence the vanishing cycles can be chosen arbitrarily. We illustrate this by introducing a family of Liouville structures on $D^6$, which depends on certain combinatorial data (given by paths on the plane, or equivalently by conjugacy classes in the braid group). There are reasons to suspect that this family should contain many exotic (nonstandard) examples, and that these are distinguished by known invariants such as symplectic cohomology. Due to the difficulty of computing these invariants, this problem (like many others of the same kind) is entirely open at present, and we present it mainly as food for thought.

\subsection{The definition}
A Liouville domain is a compact manifold with boundary $M^{2n}$, together with a one-form $\theta$ which has the following two properties. First, $\o = d\theta$ should be symplectic. Secondly, the vector field $Z$ defined by $i_Z\o = \theta$ should point strictly outwards along $\partial M$.

\begin{example} \label{th:stein}
Let $U$ be a Stein manifold, with complex structure $J$. One of the equivalent formulations of the Stein property is that there is an exhausting function $h: U \rightarrow \R$ which is (strictly) plurisubharmonic, meaning that $-dd^ch = -d(dh \circ J)$ is a K{\"a}hler form. Then, if $C$ is a regular value of $h$, the sublevel set
\begin{equation} \label{eq:sublevel}
M = h^{-1}((-\infty;C]),
\end{equation}
equipped with $\theta = -d^ch$, is a Liouville domain. The Liouville flow in this case is just the gradient flow of $h$, defined with respect to the K{\"a}hler metric associated to $\omega = d\theta$.
\end{example}

Returning to the general discussion, note that $\alpha = \theta|\partial M$ is a contact one-form on $\partial M$. The flow of $Z$ is always defined for negative times, and gives rise to a canonical collar
\begin{equation} \label{eq:collar}
\begin{aligned}
 & \kappa: (-\infty;0] \times \partial M \longrightarrow M, \\
 & \kappa^*\theta = e^r\alpha, \;\; \kappa^*Z = \partial_r,
\end{aligned}
\end{equation}
modelled on the negative half of the symplectization of $(\partial M,\alpha)$. It is therefore natural to create more space by attaching an infinite cone, which corresponds to the positive half:
\begin{equation}
\label{eq:completion}
\begin{aligned}
 & \hat{M} = M \cup_{\partial M} ([0;\infty) \times \partial M), \\
 & \hat{\theta}|([0;\infty) \times \partial M) = e^r\alpha, \;\;
 \hat{Z}|([0;\infty) \times \partial M) = \partial_r, \;\; \hat\omega = d\hat\theta.
\end{aligned}
\end{equation}
This process is called completion (because the extended Liouville field $\hat{Z}$ is complete: its flow exists for all times).

A Liouville isomorphism between domains $M_0,M_1$ is a diffeomorphism $\phi: \hat{M}_0 \rightarrow \hat{M}_1$ satisfying $\phi^*\hat\theta_1 = \hat\theta_0 + d(\text{\it some compactly supported function})$.
Obviously, any such $\phi$ is symplectic, and compatible with the Liouville flow at infinity. This means that on $[\rho;\infty) \times \partial M_0 \subset \hat{M}_0$ for some $\rho \gg 0$, it has the form
\begin{equation} \label{eq:induced-symplecto}
\phi(r,y) = (r-f(y),\psi(y)),
\end{equation}
where $\psi: \partial M_0 \rightarrow \partial M_1$ is a contact isomorphism, satisfying $\psi^*\alpha_1 = e^f \alpha_0$ for some function $f$. Note that while the contact structure at the boundary is preserved under Liouville isomorphism, the contact one-form is not, and in fact can be changed arbitrarily.
Another fact which one should keep in mind is a version of Moser's Lemma, which says that deformation equivalence implies Liouville isomorphism:

\begin{lemma} \label{th:moser}
Let $(\theta_t)_{0 \leq t \leq 1}$ be a family of Liouville structures on $M$. Then all the $(M,\theta_t)$ are mutually Liouville isomorphic.
\end{lemma}

\begin{example} \label{th:unique-sublevel}
For an $h: U \rightarrow \R$ as in the previous example, suppose that the critical point set of $h$ is compact. Then, if we take $C$ to be bigger than the largest critical value, the resulting Liouville domain \eqref{eq:sublevel} is independent of the particular choice of $C$ up to Liouville isomorphism. If one assumes in addition that $\nabla h$ is complete, $(U,-dd^ch)$ itself will be symplectically isomorphic to $\hat{M}$.
In this context, it is maybe useful to know that completeness of the gradient vector field can always be achieved by a reparametrization $h \mapsto \beta(h)$ \cite[Lemma 3.1]{biran-cieliebak01}.
\end{example}

\subsection{Lefschetz fibrations\label{subsec:lefschetz}}
We will not give a proper explanation of the theory of Lefschetz fibrations for Liouville domains, and instead just regard it as a ``black box'', which constructs symplectic manifolds from lower-dimensional data. More specifically, let $M$ be a four-dimensional Liouville domain, and $(V_1,\dots,V_m)$ an ordered collection of embedded Lagrangian two-spheres in $M \setminus \partial M$. There is a unique (up to deformation) Lefschetz fibration over the disc $D$, whose fibre is $M$ and whose vanishing cycles are the $V_k$. The total space of this fibration (after rounding off corners) is a six-dimensional Liouville domain $E$. Symplectically (up to Liouville isomorphism), $E$ depends on the Lagrangian isotopy classes of the $V_k$; but topologically (up to diffeomorphism), only the differentiable isotopy classes matter. This is important because the discrepancy between Lagrangian and differentiable isotopy is known to be large in many cases \cite{seidel98b}.

\begin{remark} \label{th:details}
This construction generalizes to Liouville domains $M$ of any dimension $2n$, except for some technical complications which we will now mention briefly. On the symplectic side, the vanishing cycles have to satisfy an exactness condition $[\theta|V_k] = 0 \in H^1(V_k;\R)$; and moreover, each such cycle should come with a diffeomorphism $f_k: S^n \rightarrow V_k$, or rather with a specified class $[f_k] \in \pi_0(\text{\it Diff}(S^n,V_k)/O_{n+1})$. The first condition is automatic unless $n = 1$, while the $\pi_0$ group appearing in the second one is trivial for $n \leq 3$.

As far as the topological side is concerned (if one wants $E$ just as a smooth manifold), what one needs is a collection of embedded $n$-spheres $V_k$ together with classes $[f_k]$ as before, and additionally isomorphisms
\begin{equation} \label{eq:nu-t}
\nu_{V_k} \iso T^*\!V_k.
\end{equation}
An easy way to see why these are relevant is to interpret the construction as attaching handles to the boundary spheres $\Sigma_k = \{z_k\} \times V_k \subset \partial(D^2 \times M)$, for some choice of distinct cyclically ordered points $z_k$. The isomorphisms \eqref{eq:nu-t} then identify the normal bundles of those boundary spheres with $\R \oplus T^*\!S^n$, which is canonically trivial. This is the framing used to attach the handles (of course, if the $\Sigma_k$ are Lagrangian, they come with distinguished maps \eqref{eq:nu-t}, hence no additional framing data is required; see Section \ref{subsec:handles} for more discussion along those lines). Once one has taken orientations into account, the remaining freedom in choosing \eqref{eq:nu-t} is an element of $\pi_n(SO_n)$. Again, this vanishes for $n \leq 2$, which is why the issue did not appear in our original discussion.
\end{remark}

To proceed to a concrete example, let's fix a polynomial $p \in \C[z]$ which is monic of degree $m+1 \geq 3$, and has no multiple zeros. Denote the set of zeros by $P \subset \C$. Consider the smooth affine algebraic surface $U \subset \C^3$ defined by
\begin{equation} \label{eq:ale}
xy = p(z),
\end{equation}
and equip it with the exact symplectic form inherited from (the standard form on) $\C^3$. To any embedded path $c \subset \C$ whose endpoints lie in $P$, and which otherwise avoids $P$, one can associate a Lagrangian two-sphere $L_c \subset X$, defined explicitly as $L_c = \bigcup_{z \in c} C_z$, where $C_z = \{(x,y,z) \in U, |x| = |y|\}$. It is an elementary obser\-vation that the differentiable isotopy class of $L_c$ depends only on the endpoints of $c$. In contrast, Floer cohomology computations from \cite{khovanov-seidel98} show that the Lagrangian isotopy class of $L_c$ recovers the isotopy class of $c$ (among paths of this kind; which means that isotopies may not pass through $P$). Define a chain of paths to be a collection $(c_1,\dots,c_m)$, with the property that for some ordering $P = \{z_1,\dots,z_{m+1}\}$, the endpoints of $c_k$ are $\{z_k,z_{k+1}\}$. Call a chain standard if the intersections are the least possible, which means that $c_k \cap c_l = \emptyset$ for $|k-l| \geq 2$, and $|c_k \cap c_{k+1}| = 1$.

\begin{figure}
\begin{centering}
\includegraphics{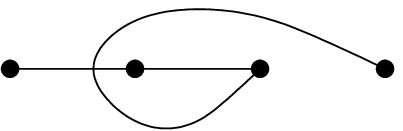}
\caption{\label{fig:chain}}
\end{centering}
\end{figure}

Fix a chain, and let $M \subset U$ be the intersection of $U$ with a large ball in $\C^3$, containing all our Lagrangian spheres $L_{c_k}$. Take those as vanishing cycles, $V_k = L_{c_k}$, and construct the resulting $E$. If the chain is standard, then $E$ is actually Liouville isomorphic to standard $D^6$ (and therefore, its completion is a standard symplectic $\R^6$). This is not obvious from the point of view taken here, but actually the resulting Lefschetz fibration is well-known (it is the Morsification of the function $xy + z^{m+1}$, known in singularity theory as the $A_m$ type singularity). On the other hand, it is true for any chain that $E$ is diffeomorphic to $D^6$. This follows from our previous discussion of isotopy, which implies that the diffeomorphism type of $E$ is the same for all chains. Alternatively, one can use the fact that the $L_{c_k}$ form a basis of $H_2(M)$ to show that $E$ is contractible, and then appeal to the h-cobordism theorem. Finally, we should point out that even non-standard chains may give rise to standard Liouville structures. Examples can be easily constructed using known moves (analogues of handle slide and handle cancellation), see for instance Figure \ref{fig:chain}. In spite of that, we propose:

\begin{conjecture}
Consider the set of chains up to homotopy, with some bound $\Lambda$ on the complexity of the paths (this can be done by taking lengths with respect to a metric, or by counting the number of pieces with respect to a suitable decomposition of the disc; the resulting set is finite). Then, as $\Lambda \rightarrow \infty$, the probability that a randomly chosen chain will give rise to a nonstandard Liouville structure on $D^6$ (and actually to a nonstandard symplectic structure on the completion $\R^6$) approaches $1$.
\end{conjecture}

\section{Symplectic cohomology}

Symplectic cohomology is Hamiltonian Floer cohomology for a specific class of (Hamiltonian) functions, which are tailored to completions of Liouville domains. Cieliebak, Floer and Hofer \cite{floer-hofer94,cieliebak-floer-hofer95} were the first ones to to introduce such ideas into pseudo-holomorphic curve theory (this drew on the insights obtained in previous work, which used a more conventional variational approach). Their constructions were made with a view to applications in ``quantitative'' symplectic topology, where the theory has indeed been highly successful; see for instance \cite{cieliebak-floer-hofer-wysocki96}. Subsequently, Viterbo \cite{viterbo97a} introduced a more ``qualitative'' version (which had an additional source of inspiration, namely results coming from the theory of generating functions). This is the definition we use here.

Since the Reeb flow on the boundary appears prominently in the definition, it is natural to start by looking at cases where that flow is reasonably well understood. A highly nontrivial manifestation of this idea is Viterbo's theorem on cotangent bundles. Besides stating this result, we also describe a few other computations, of a somewhat simpler nature. We then turn to invariance under Liouville isomorphism, which is the most important basic property of symplectic cohomology. The proof outlined here is quite typical of the general style of the theory: it proceeds via a sequence of ``approximate isomorphisms'', constructed using continuation maps, which finally ``converge to'' an actual isomorphism when passing to the direct limit. Will see in the next section that this direct limit argument is not only a useful technical device, but also gives rise to additional invariants.

\subsection{A first look}
Fix an abelian group $\K$. Let $M$ be a Liouville domain of dimension $2n$. The symplectic cohomology of $M$ with $\K$-coefficients is a $\Z/2$-graded $\K$-vector space $SH^*(M)$, which comes with a natural $\Z/2$-graded map
\begin{equation} \label{eq:classical}
H^{*+n}(M;\K) \longrightarrow SH^*(M).
\end{equation}
This map is the central object of interest in \cite{viterbo97a,viterbo97b}, because its failure to be an isomorphism signals the presence of at least one periodic orbit of the Reeb flow on $\partial M$. Hence, it can be used to study Weinstein's conjecture on the existence of such orbits (this is an efficient line of attack, but it has its limits, since it depends on finding suitable Liouville fillings of a given contact manifold).

Often, it is convenient to upgrade the $\Z/2$-grading to a $\Z$-grading. For that, assume that $c_1(M) = 0$. More precisely, after choosing a compatible almost complex structure, $-c_1(M)$ is represented canonically by a complex line bundle $\scrK$, and we want to choose a trivialization of that line bundle, up to homotopy. This then gives rise to a $\Z$-grading of $SH^*(M)$, which is such that \eqref{eq:classical} is of degree zero. Of course, if $H^1(M) = 0$, there is only one homotopy class of trivializations, and then no choice is necessary. Otherwise, the grading may depend on the trivialization (actually, it does so in a precisely specified way, depending on the splitting of $SH^*(M)$ into pieces parametrized by free homotopy classes of loops in $M$).

\begin{remark}
It is actually possible to define symplectic cohomology for more general symplectic manifolds with contact type boundary (ones where the Liouville flow does not extend over the interior). As usual, one either has to make some assumptions governing holomorphic spheres, or else work with $\Q$-coefficients. This generalization is outside the scope of our discussion.
\end{remark}

One of the most elementary cases is when the Reeb flow on $\partial M$ is a circle action (every orbit is periodic). To make things even simpler, let's assume that $\scrK$ is trivial, that the circle action is free, and that $\partial M$ is connected. Then there is a spectral sequence converging to (the $\Z$-graded version of) $SH^*(M)$, whose starting term is
\begin{equation} \label{eq:morse-bott}
E_1^{pq} = \begin{cases}
H^{q+n}(M;\K) & p = 0, \\
H^{p+q+n-p\mu}(\partial M;\K) & p< 0, \\
0 & p>0.
\end{cases}
\end{equation}
The map \eqref{eq:classical} appears naturally in this context, as the edge homomorphism of the spectral sequence. The integer $\mu \in 2\Z$ is a Conley-Zehnder type index, which can be explained as follows. Take a primitive orbit $x$ of our circle action. The trivialization of $\scrK$ determines a symplectic trivialization of $x^*TM$ (up to homotopy). On the other hand, the circle action itself (extended to a collar neighbourhood of $\partial M$ in the obvious way) gives a family of isomorphisms $TM_{x(0)} \rightarrow TM_{x(t)}$. In view of our trivialization, these turn into a loop in $Sp_{2n}$, and $\mu/2$ is the class of this loop in $\pi_1(Sp_{2n}) \iso \Z$. For future reference, we point out that in the simplest case $H^1(\partial M) = 0$, the choice of trivialization of $\scrK$ is irrelevant, so that $\mu$ only depends on $(\partial M,\alpha)$. The construction of the spectral sequence is relatively straightforward starting from the definition of $SH^*(M)$, requiring no more than standard Bott-Morse methods \cite{pozniak,bourgeois02}.

\begin{example} \label{th:ball}
Let $M = B^{2n}$ be the unit ball in $\R^{2n}$. This satisfies $\mu = 2n$, hence the $E_1$ term of the spectral sequence has one generator in each degree $-n,-n-1,-3n,-3n-1,-5n,\cdots$. It turns out that the differential $d_1$ on this page is acyclic, so that $SH^*(M) = 0$. We will see several proofs of this vanishing result later on.
\end{example}

\begin{example} \label{th:surface}
Take $M$ to be a surface of genus $g>0$ with one boundary component. Then $\mu = 2\chi(M) = 2-4g < 0$. For topological reasons (the boundary and its multiple covers are non-contractible), the differentials in the spectral sequence must vanish. Hence, $SH^*(M)$ is the direct sum of $H^*(M;\K)$ shifted down by $1$, and copies of $H^*(S^1;\K)$ shifted up by $4g-3, 8g-5, \cdots$.
\end{example}

\begin{example} \label{th:cotangent-sphere}
Another example is $M = D^*S^n$, the cotangent disc bundle of the $n$-sphere, formed with respect to the round metric (cotangent bundles of oriented manifolds $N$ always come with a preferred trivialization of $\scrK$; topologically speaking, the reason is that the structure group of $T(T^*N)$ is $SO_n \subset Sp_{2n}$). Then $\mu = 2n-2$, so the $E_1$ term consists of a copy of $H^*(S^n;\K)$ shifted down by $n$, plus copies of $H^*(\partial D^*S^n;\K)$ (the unit sphere tangent bundle) shifted down by $3n-2, 5n-4, \dots$. Let's take $\K = \Q$ for simplicity. Because of the gradings, and of the direction in which the spectral sequence differentials point, it is clear that the spectral sequence must degenerate for $n > 3$.
\end{example}

Actually, in the last-mentioned example it is true that the spectral sequence degenerates for all $n$. This cannot be derived from mere grading considerations, and it follows instead from a fundamental theorem of Viterbo \cite{viterbo97b}, which determines the symplectic cohomology of general cotangent bundles. Besides its considerable intrinsic importance, this theorem has received a lot of attention recently, because it establishes a connection (partly proved, and partly still conjectural) between symplectic cohomology and string topology.
The statement is:

\begin{theorem}[Viterbo] \label{th:viterbo}
Let $M = D^*N$ be the cotangent disc bundle of an oriented closed manifold. Then $SH^*(M) \iso H_{-*}({\mathcal L} N;\K)$ is the homology of the free loop space, with the grading reversed. Under this isomorphism, the map \eqref{eq:classical} turns into the inclusion of constant loops,
\begin{equation}
H^{n+*}(M;\K) = H^{n+*}(N;\K) \iso H_{-*}(N;\K) \longrightarrow H_{-*}({\mathcal L} N;\K).
\end{equation}
\end{theorem}

Intuitively, the basic idea is that the Reeb flow on $\partial M$ agrees with the geodesic flow on $N$ (for the Riemannian metric used to define the disc bundle inside $T^*N$). Periodic Reeb orbits correspond to closed geodesics, and in view of the classical Morse theory for the geodesic energy functional, one is naturally led to look at the homology of the free loop space. There are a few simple cases where the theorem can be proved in an elementary way just by comparing the generators in the two Morse-type complexes; for instance, if $N$ admits a Riemannian metric with negative sectional curvature. However, the general statement of Theorem \ref{th:viterbo} is beyond the reach of such methods; for instance, it is by no means obvious why $SH^*(D^*N)$ should be a homotopy invariant of $N$. Besides the original proof \cite{viterbo97b}, there are now two other (analytically more complicated, but overall more direct) approaches, by Abbondandolo-Schwarz \cite{abbondandolo-schwarz06} and Salamon-Weber \cite{salamon-weber03}, respectively.

\subsection{Floer theory basics}
Let $M$ be any symplectic manifold whose symplectic form is exact, $\omega = d\theta$. Fix $H \in \smooth(M,\R)$, and let $X$ be the associated Hamiltonian vector field, which in our convention satisfies $\omega(\cdot,X) = dH$. The $H$-perturbed action functional on the free loop space ${\mathcal L} M$ is the function
\begin{equation} \label{eq:action}
 A_H(x) = \,-\!\!\int x^*\theta + \int H(x(t))\,dt.
\end{equation}
Critical points of $A_H$ are precisely those $x$ which are flow lines of $X$, which means they are $1$-periodic orbits of that vector field. Now let $J$ be an $\omega$-compatible almost complex structure. Floer's equation (in its simplest) form is the PDE
\begin{equation} \label{eq:floer}
\left\{
\begin{aligned}
& u: \R \times S^1 \longrightarrow M, \\
& \partial_s u + J(\partial_t u - X) = 0.
\end{aligned}
\right.
\end{equation}
The energy of a solution is, by definition,
\begin{equation} \label{eq:energy}
E(u) = \int |\partial_s u|^2 \, ds \wedge dt = \int \omega(\partial_s u,\partial_t u - X)\, ds \wedge dt.
\end{equation}
Recall that \eqref{eq:floer} is, formally, the negative gradient flow equation for $A_H$. It is therefore natural to look at connecting trajectories, which means solutions with asymptotic behaviour $\lim_{s \rightarrow \pm\infty} u(s,\cdot) = x_{\pm}(\cdot)$, where $x_{\pm}$ are $1$-periodic orbits of $X$. For any such trajectory, we have the a priori energy estimate
\begin{equation} \label{eq:a-priori}
E(u) = A_H(x_-) - A_H(x_+).
\end{equation}

A common variation on \eqref{eq:floer} is the continuation map equation, used to construct homomorphisms between Floer cohomology groups. This involves families $\{H_s\}$ and $\{J_s\}$ of functions and almost complex structures, depending on $s \in \R$, which should be eventually constant: $(H_s,J_s) = (H_-,J_-)$ for $s \ll 0$, $(H_s,J_s) = (H_+,J_+)$ for $s \gg 0$. The equation is
\begin{equation} \label{eq:cont}
\partial_s u + J_s(\partial_t u - X_s) = 0,
\end{equation}
and the natural asymptotic condition is $\lim_{s \rightarrow \pm \infty} u(s,\cdot) = x_{\pm}(\cdot)$, where $x_{\pm}$ is a $1$-periodic orbit of $X_{\pm}$. If one defines the energy as in \eqref{eq:energy}, then a solution with such limits satisfies
\begin{equation} \label{eq:cont-energy}
 E(u) = A_{H_-}(x_-) - A_{H_+}(x_+) + \int (\partial_s H_s)(u) \, ds \wedge dt.
\end{equation}
For instance, if $\partial_s H_s \leq 0$ everywhere, the a priori bound will be as good as (in fact better than) \eqref{eq:a-priori}. On the other hand, there are situations where $\partial_s H_s$ is not bounded above, leading to a failure of compactness, which means that continuation maps can not be defined. This is a fundamental point, even though in our subsequent discussion, it will be somewhat obscured by the more technical device of maximum principles.

\subsection{Symplectizations\label{subsec:symplectizations}}
Getting somewhat closer to the intended applications, let's consider the case where the target $M = \R \times Y$ is the symplectization of a contact manifold $(Y,\alpha)$. We will be interested in Hamiltonian functions of the form
\begin{equation} \label{eq:H-h}
H = h(e^r),
\end{equation}
where $r$ is the first variable on $M$, and $h$ any smooth function with $h' \geq 0$. By definition, the Hamiltonian vector field of $e^r$ is the Reeb field $R$ of $\alpha$ (pulled back to $M$ in the obvious way). Hence, the vector field of a general function \eqref{eq:H-h} is $X = h'(e^r)R$. Clearly, there is a close relation between the periodic orbits of $X$ and $R$. Specifically, if $y(t)$ is a $T$-periodic Reeb orbit in $Y$, and $r \in \R$ a number satisfying $h'(e^r) = T$, then
\begin{equation} \label{eq:x-orbit}
x(t) = (r,y(Tt))
\end{equation}
is a $1$-periodic orbit of $X$. These are all $1$-periodic orbits of $X$ except for the stationary ones, which occur where $h'(e^r) = 0$. By inserting this into \eqref{eq:action}, one sees that the value of $A_H$ at a critical point is
\begin{equation} \label{eq:reeb-action}
 A_H(x) = h(e^r) - e^r h'(e^r).
\end{equation}
For example, suppose that we take a function $h$ satisfying $\lim_{r \rightarrow -\infty} h'(e^r) = 0$, $\lim_{r \rightarrow +\infty} h'(e^r) = \infty$, and $h''(e^r) > 0$ everywhere. In that case, periodic Reeb orbits of any period correspond bijectively to $1$-periodic orbits of $X$ (note that the orbits are considered to be parametrized, and that multiply-covered orbits are also allowed). More precisely, Reeb orbits with larger periods $T$ correspond to $1$-periodic orbits \eqref{eq:x-orbit} with bigger values of $r$, and smaller values of the action functional (since by assumption, $\partial_r(h(e^r) - e^r h'(e^r)) = -e^{2r} h''(e^r) < 0$).

The non-compactness of the target manifold gives rise to some additional issues. Specifically, we will often need to show that sequences of solutions $u$ do not escape to infinity in the $r \rightarrow +\infty$ direction (unlike the case of contact homology, the other end $r \rightarrow -\infty$ is irrelevant, because it will eventually be capped off by a compact Liouville domain). Assume that the almost complex structure $J$ is of contact type, which means that
\begin{equation} \label{eq:contact-type}
d(e^r) \circ J = -\theta.
\end{equation}
Equivalently, in terms of the $\omega$-orthogonal splitting $T(\R \times Y) = (\R \partial_r \oplus \R R) \oplus \xi \iso \C \oplus \xi$, these $J$ are standard on the first summand (and arbitrary, except for the $d\alpha$-compatibility condition, on the second summand; note that this second component is allowed to vary with $r$). Let $u$ be a solution of Floer's equation, and consider the function $\rho = e^r \circ u: \R \times S^1 \rightarrow \R$. In view of \eqref{eq:contact-type}, this satisfies
\begin{align}
&
\begin{aligned}
 & \partial_s \rho = \theta(\partial_t u) - \rho h'(\rho), \\
 & \partial_t \rho = -\theta(\partial_s u), \\
\end{aligned} \\
\intertext{which can be written in a less coordinate-bound way as}
 \label{eq:dc-term} & d^c\rho = d\rho \circ i = -u^*\theta + \rho\cdot h'(\rho)\, dt. \\
\intertext{By differentiating again, and substituting $|\partial_s u|^2 = \omega(\partial_s u,\partial_t u - X) = \omega(\partial_s u,\partial_t u) - dH(\partial_s u) =  \omega(\partial_s u,\partial_t u) - h'(\rho) \cdot \partial_s \rho$, one gets}
 & \Delta \rho = |\partial_s u|^2 - \rho\cdot h''(\rho)\cdot \partial_s \rho. \label{eq:maximum}
\end{align}
Like any solution of an inequality $\Delta \rho + v(s,t) \cdot \partial_s \rho + w(s,t) \cdot \partial_t \rho \geq 0$, our function $\rho$ obeys a maximum principle. This means that for any bounded open subset $\Omega \subset \R \times S^1$, the maximum of $u|\bar\Omega$ must occur on the boundary. To give a sample application, suppose that $u$ is a a connecting trajectory between $x_- = (r_-,y_-)$ and $x_+ = (r_+,y_+)$. Then, the entire image of $u$ must be contained in the subset where $r \leq \max(r_-,r_+)$.

Analogously, one can consider solutions of \eqref{eq:cont}, where $H_s(r,y) = h_s(e^r)$, and each $J_s$ is of contact type at infinity. If one defines $\rho$ as before, the computation goes through in the same way until \eqref{eq:dc-term}, but the final differentiation creates an additional term
\begin{equation} \label{eq:cont-delta}
 \Delta \rho = |\partial_s u|^2 - \rho\cdot h''_s(\rho)\cdot \partial_s \rho - \rho\cdot (\partial_s h'_s)(\rho).
\end{equation}
Bearing in mind that $\rho$ is positive by definition, one finds that the maximum principle only applies if $\partial_s h'_s \leq 0$. This condition is similar to (but not the same as) the one we encountered before, when trying to get a priori bounds from \eqref{eq:cont-energy}.

Occasionally, it will be useful to extend the class of functions and almost complex structures under consideration, in order to make it invariant under symplectic isomorphisms \eqref{eq:induced-symplecto}. Namely, one sets $R(r,y) = r-f(y)$ for some $f \in \smooth(Y,\R)$, and instead of \eqref{eq:H-h}, \eqref{eq:contact-type}, requires that
\begin{equation} \label{eq:changed}
\begin{aligned}
 & H = h(e^R), \\
 & d(e^R) \circ J = -\theta.
\end{aligned}
\end{equation}
As far as Floer's equation is concerned, this is not actually more general than the previous framework, to which it reduces after changing coordinates and the contact one-form: $(r,y) \mapsto (R(r,y),y)$, $\alpha \mapsto e^f\alpha$. However, in the context of continuation maps, one can now take $f_s$ (and hence $R_s$) which varies with $s$, assuming as usual that it is locally constant for $|s| \gg 0$. Correspondingly, take $H_s$ and $J_s$ to be of the form \eqref{eq:changed}. Given a solution $u$ of the resulting equation \eqref{eq:cont}, one considers
\begin{equation} \label{eq:strange-rho}
\rho(s,t) = \exp(R_s(u(s,t))).
\end{equation}
The same kind of computation as before yields
\begin{equation} \label{eq:delta-rho-big}
\begin{aligned}
\Delta \rho  = \;&
|\partial_s u|^2 - \rho \cdot h_s''(\rho) \cdot \partial_s \rho - \partial_s f_s \cdot \partial_s \rho \\
 & - \rho \cdot (\partial_s h_s')(\rho)
 + \rho \cdot h_s'(\rho) \cdot \partial_s f_s - \rho \cdot \partial_s^2 f_s
 - \rho \cdot d(\partial_s f_s)(\partial_s u).
\end{aligned}
\end{equation}
It is implicit in our notation that $f_s$ and its derivatives are always evaluated at $u(s,t)$. In applications, the family $\{f_s\}$ is given, and one wants to choose $\{h_s\}$ in such a way that the maximum principle applies. The most ``dangerous'' term in \eqref{eq:delta-rho-big} is the last one, which is potentially an unbounded multiple of $\rho$. However, thanks to the exponential growth of the metric on $M$, $|\rho \cdot d(\partial_s f_s)(\partial_s u)|$ is actually bounded above by $C|\partial_s u|$. Here and later, $C$ stands for some large constant depending only on $\{f_s\}$ (which can be a different one each time the notation occurs). By exploiting this, one gets an inequality of the form
\begin{equation} \label{eq:strange-max}
\begin{aligned}
 \Delta \rho + (\rho \cdot h_s''(\rho) + \partial_s f_s)\cdot \partial_s\rho & \geq
 \rho \big(-\partial_s h_s'(\rho) - C  h_s'(\rho) - C) - C\,|\partial_s u| + |\partial_s u|^2 \\
 & \geq \rho \big(-\partial_s h_s'(\rho) - C h_s'(\rho) - C) - C.
\end{aligned}
\end{equation}
Wanting the right hand side to be positive, at least on the subset where $\rho$ is large, roughly comes down to an exponential decay condition on $h_s'$. Of course, this only needs to hold on the bounded subset of those $s$ where $\{f_s\}$ is not constant; elsewhere, the behaviour of $\rho$ is governed by the equation \eqref{eq:cont-delta} for ordinary continuation maps, so that $\partial_s h_s' \leq 0$ suffices.

\subsection{First definition}
Let $\hat{M}$ be the completion of a Liouville domain $M$. Choose some compatible almost complex structure $J$ on this, which is of contact type at infinity. Similarly, we will consider Hamiltonian function which at infinity are of the form $H(r,y) = h(e^r)$, where the function $h$ satisfies $\lim_{r \rightarrow \infty} h'(e^r) = \infty$. By definition, symplectic cohomology is the associated Floer cohomology:
\begin{equation} \label{eq:sh-1}
SH^*(M) = HF^*(H).
\end{equation}
Roughly speaking, this means that it is the Morse cohomology of the action functional \eqref{eq:action}, so the underlying chain complex $CF^*(H)$ should be generated by critical points, and the differential $\delta$ given by counting solutions of Floer's equation. To explain the conventions here, note first that by ``generated'' we mean that $CF^*(H)$ is constructed as a direct sum (not direct product); in finite-dimensional terms, the model is cohomology with compact supports. Secondly, since we are talking about cohomology, and \eqref{eq:floer} is the negative gradient flow equation, a solution with asymptotics $x_{\pm}$ contributes a term to $\delta$ which maps $x_+ \mapsto x_-$, up to sign.

On the technical side, note that the critical points of the action functional are usually degenerate, because of the rotational symmetry (recall that we are considering parametrized periodic orbits of the Reeb flow, which obviously occur in $S^1$ families). The standard way to resolve this difficulty is to break the symmetry by making a small $t$-dependent perturbation of $H$. Similarly, transversality issues for the moduli spaces of connecting trajectories are addressed by taking a $t$-dependent $J$. To simplify the presentation, we will systematically suppress any mention of such perturbations, and just pretend that we can work with the degenerate data itself. This is of course incorrect, but is hopefully permitted in an informal exposition such as the present one.

With this in mind, the maximum principle derived from \eqref{eq:maximum} implies that given endpoints $x_\pm$, all solutions of \eqref{eq:floer} remain within a compact subset of $\hat{M}$. This neutralizes all issues arising from the non-compactness of the target space, ensuring that our Floer complex is well-defined. To show that its cohomology is independent of the choice of $(H,J)$, one can look at a one-parameter family connecting two choices, and then study the bifurcations occurring in the Floer complex, in the manner of \cite{floer88c,lee01}. This kind of argument is fairly direct, but the underlying analytic work is substantial (and it seems even harder to show that the isomorphisms are canonical).

\subsection{Second definition\label{subsec:limit}}
Most of the literature actually favours an alternative approach, based on direct limits. Fix some $\tau>0$ which is not the period of any Reeb orbit on $\partial M$ (this is a generic condition, since the periods form a countable closed subset of $\R^+$), and take a Hamiltonian function $H^\tau$ which at infinity satisfies $H^\tau(r,y) = \tau e^r + \text{\it (constant)}$. Define
\begin{equation}
SH^*(M)^{<\tau} = HF^*(H^\tau).
\end{equation}
The notation reflects the fact that this group takes into account only Reeb orbits of period $<\tau$ (one can extend the definition to general $\tau$, by taking a function whose slope approaches $\tau$ from below, but we will not really need this). Now suppose that we have two such values $\tau_{\pm}$, and corresponding functions $H_{\pm} = H^{\tau_{\pm}}$. Choose a family $\{\tau_s\}$ interpolating between $\tau_s = \tau_-$ for $s \ll 0$ and $\tau_s = \tau_+$ for $s \gg 0$. Take a corresponding family $\{H_s = H^{\tau_s}\}$ of functions, as well as a family $\{J_s\}$ of almost complex structures, and consider the continuation map equation \eqref{eq:cont}. In view of \eqref{eq:cont-delta}, the maximum principle for solutions will hold outside a compact subset, provided that
\begin{equation} \label{eq:tau-decreases}
\partial_s \tau_s \leq 0.
\end{equation}
Note that once the maximum principle is valid, a bound on $\int (\partial_s H_s)(u) \, ds \wedge dt$ follows immediately because $u$ must remain within a compact subset; this means that controlling the energy \eqref{eq:cont-energy} is unproblematic. The first application of the resulting continuation maps, in the special case where $\tau_+ = \tau_-$, is to give an elegant proof that $SH^*(M)^{<\tau}$ is independent of all choices, up to canonical isomorphism; this follows the standard Floer theory strategy. More interestingly, one gets canonical maps
\begin{equation}
SH^*(M)^{<\tau_+} \longrightarrow SH^*(M)^{<\tau_-} \quad \text{for all $\tau_+ < \tau_-$.}
\end{equation}
These behave well under composition, hence form a direct system indexed by all possible $\tau$. The second definition of symplectic cohomology is as the limit of this system:
\begin{equation} \label{eq:sh-2}
SH^*(M) = \underrightarrow{lim}_\tau\, SH^*(M)^{<\tau} = \underrightarrow{lim}_\tau\, HF^*(H^\tau).
\end{equation}
To see that this agrees with \eqref{eq:sh-1}, one again uses continuation maps, but now for families of functions $\{H_s\}$ which satisfy $H_s(r,y) = h_s(e^r)$ at infinity, where at one end $s \ll 0$, $h_s = h$ has unbounded growth $\lim_{s \rightarrow} h'(r) = \infty$, and at the other end $s \gg 0$, $h_s(e^r) = \tau e^r + \text{\it (constant)}$. This gives rise to a chain maps
\begin{equation} \label{eq:infinite-continuation}
CF^*(H^\tau) \longrightarrow CF^*(H)
\end{equation}
which, in the limit $\tau \rightarrow \infty$, induce a map between from the right hand side of \eqref{eq:sh-2} to \eqref{eq:sh-1}. A suitable choice of functions ensures that (for a sequence of $\tau$'s going to $\infty$) $CF^*(H^{\tau})$ is a subcomplex of $CF^*(H)$, and that the map \eqref{eq:infinite-continuation} is the inclusion of that subcomplex. As $\tau \rightarrow \infty$, these complexes exhaust $CF^*(H)$, hence we get an isomorphism (here, one sees why it is important to define $CF^*(H)$ as a direct sum).

It is instructive to see how the map \eqref{eq:classical} appears in either context. If one takes \eqref{eq:sh-2} as the definition, taking $\tau$ to be very small (at least, smaller than the length of the shortest periodic Reeb orbit) yields a group $HF^*(H^\tau)$ which, by standard Floer-theoretic arguments, is canonically isomorphic to the ordinary cohomology $H^{*+n}(M;\K)$. Note that the absolute cohomology of $M$, rather than the relative cohomology of $(M,\partial M)$, appears, because the function $H^\tau$ grows at infinity. In these terms, \eqref{eq:classical} is simply one of the maps $HF^*(M)^{<\tau} \rightarrow HF^*(M)$ which are obviously part of the direct limit formalism. If one takes \eqref{eq:sh-1} as a starting point, the argument is similar but a little more delicate. Suppose that we set $H = 0$ on $M$, and $H(r,y) = h(e^r)$ on the entire cone $[0;\infty) \times \partial M$, where $h'(e^r) > 0$ for all $r>0$. Then, the stationary orbits of  $X$ corresponding to points of $M$ have action $A_H = 0$, whereas the $1$-periodic orbits coming from periodic orbits of $R$ have strictly negative action $A_H<0$, by \eqref{eq:reeb-action}. Even after a small perturbation which makes the critical points nondegenerate, one still has a distinguished subcomplex of $CF^*(H)$ corresponding to generators whose action is $>-\epsilon$, for some small $\epsilon$. Again applying standard methods from Floer theory, one can show that for a suitable choice of perturbation, this subcomplex can be identified with a Morse complex computing $H^{*+n}(M;\K)$.

Finally, suppose that we have two isomorphic Liouville domains $M_{\pm}$, maximal periods $\tau_{\pm}$, and associated groups $SH^*(M_\pm)^{< \tau_\pm}$. After identifying their completions, this is the same as having a single $M$, but considering functions $H_\pm$ on $\hat{M}$ which at infinity are of the form $H_\pm(r,y) = \tau_{\pm} e^{r-f_\pm} + \text{\it (constants)}$, for some $f_\pm \in \smooth(\partial M,\R)$ (of course, by making the identification in an appropriate way, one can get either $f_+$ or $f_-$ to vanish, but we prefer the more general notation for its greater symmetry). This puts us in the situation previously considered in \eqref{eq:changed}. To define continuation maps, one needs to take a family $\{f_s\}$ interpolating between $f_\pm$, and chooses functions
\begin{equation} \label{eq:tau-s}
H_s(r,y) = \tau_s e^{r-f_s} + \text{\it (constants)}
\end{equation}
so that the maximum principle in \eqref{eq:strange-max} works out. For any given $\tau_+$, this can be done provided that $\tau_- \gg \tau_+$. In the limit \eqref{eq:sh-2}, one gets a map $SH^*(M_+) \rightarrow SH^*(M_-)$, and similarly an inverse isomorphism.

\subsection{The ball}
We now return to the case of $M = D^{2n}$, where (as stated without proof in Example \ref{th:ball}) symplectic cohomology vanishes. One can identify $\hat{M}$ with $\C^{n}$, equipped with the standard one-form $-d^c(\quarter |x|^2)$. To define $SH^*(M)$ as in \eqref{eq:sh-1}, take any function of the form
\begin{equation}
H(x) = h(\half |x|^2),
\end{equation}
where $\lim_{r \rightarrow \infty} h'(r) = \infty$. Suppose that $h'(0) \in (2\pi k;2\pi (k+1))$ for some integer $k \geq 0$, and that $h'' > 0$ everywhere. In that case, $X$ has a stationary point at $x = 0$, and all the other $1$-periodic orbits occur in the spheres $\half |x|^2 = r$, where $h'(r) \in 2\pi l$ for $l > k$. The Conley-Zehnder indices of these orbits, including the stationary one, depend on $k$. More precisely, after a small perturbation, we get a chain complex $CF^*(H)$ which has one generator in each of the following degrees:
\begin{equation}
-n-2nk,-n-2nk-1,-n-2n(k+1),-n-2n(k+1)-1,-n-2n(k+2),\cdots
\end{equation}
Clearly, by taking $k$ sufficiently large, one can get $SH^*(M) = HF^*(H)$ to be zero in any given degree. Because the group is independent of the particular choice of Hamiltonian, it follows that it must vanish altogether (readers interested in understanding this geometrically might want to think about the cancellation which occurs when $h'(0)$ crosses one of the values in $2\pi\Z$).

We should also mention a variant of this argument, based on \eqref{eq:sh-2}. Namely, for each $k$ pick some $\tau_k \in (2\pi k;2\pi (k+1))$, and consider the function
\begin{equation}
H^{\tau_k}(x) = \tau_k \cdot \half|x|^2.
\end{equation}
Clearly, $HF^*(H^{\tau_k})$ is one-dimensional, with the single generator corresponding to the stationary point $x = 0$. Moreover, one can clearly find a family $\{H_s\}$ interpolating between $H_- = H^{\tau_{k+1}}$ and $H_+ = H^{\tau_k}$, which satisfies $\partial_s H_s \leq 0$ everywhere. In this case, the energy inequality \eqref{eq:cont-energy} for solutions $u$ of the continuation map equation yields $E(u) \leq 0$, which means that the only possibility is the constant solution $u(s,t) \equiv 0$. An index computation shows that this solution is not regular: it belongs to a moduli space with virtual dimension $-2n < 0$, and will therefore disappear after a generic small perturbation, leaving an empty moduli space of solutions. As a consequence, the maps $SH^*(M)^{<\tau_k} \rightarrow SH^*(M)^{<\tau_{k+1}}$ vanish, which yields zero as the direct limit. Note that this uses only the local virtual dimension for solutions of the continuation map equation (which always makes sense), rather than the degrees of generators in the Floer complex (which only make sense when $c_1 = 0$). Of course, in the case of $M = D^{2n}$ this is merely a cosmetic advantage, but it becomes relevant in generalizations.

\section{Growth measures and affine varieties\label{sec:affine}}

Liouville domains are most obviously related to complex analytic (Stein) manifolds. Still, one can ask whether it makes any difference if our manifold is (affine) algebraic. Intuitively, one answer would be that in the algebraic case, the geometry of the boundary is ``tame'', due to the existence of a smooth compactification with normal crossings at infinity.

On the symplectic cohomology side, this connects with the following train of thought. Even in the $\Z$-graded framework, it happens often that each group $SH^k(M)$ is infinite-dimensional, and in that case, its information content as an invariant is rather limited. One can try to combat this problem by considering additional algebraic structures (discussed briefly in Section \ref{sec:operations}), but those tend to be hard to compute. In contrast, each group $SH^*(M)^{<\tau}$ is clearly finite-dimensional. For fixed $\tau$, these groups are not Liouville invariants, but by looking at the entire direct system, one can extract some quantitative information, which measures the degree of (polynomial) growth.

\subsection{Growth rates\label{subsec:growth}}
We begin by revisiting the framework set up at the end of Section \ref{subsec:limit} for proving Liouville invariance. Namely, we are given a family $\{f_s\}$ on $\partial M$ interpolating between $f_\pm$, and want to choose $\{H_s\}$ as in \eqref{eq:tau-s} so that solutions of the associated continuation map equation remain within a bounded subset. This means that the maximum principle has to apply to \eqref{eq:strange-rho}, at least in the region where $\rho$ is large. In view of \eqref{eq:strange-max}, this means that
\begin{equation} \label{eq:decay}
 \begin{cases}
 \partial_s \tau_s \leq -C \tau_s - C & \text{on a bounded interval $I \subset \R$}, \\
 \partial_s \tau_s \leq 0 & \text{elsewhere.}
\end{cases}
\end{equation}
The constant $C$ and the interval $I$ depend only on $\{f_s\}$. Integrating out \eqref{eq:decay} yields an inequality of the form $\tau_- + 1 \geq e^{\mathrm{length}(I) C} (\tau_+ + 1)$. To simplify this slightly, take some $D > e^{\mathrm{length}(I) C}$, so that at least for $\tau_+ \gg 0$, the condition $\tau_- \geq D \tau_+$ is sufficient. Going back to the original problem, this means that we have continuation maps $SH^*(M_+)^{<\tau} \rightarrow SH^*(M_-)^{<D \tau}$ for all $\tau \gg 0$. It is easy to check that these form a map of direct systems. The same argument applies if we exchange $M_+$ and $M_-$, yielding a ladder-shaped commutative diagram
\begin{equation} \label{eq:ladder}
\xymatrix{
\cdots & \cdots \\
\ar[u]
SH^*(M_+)^{<D^4\tau} \ar[r] & SH^*(M_-)^{<D^5\tau} \ar[ul] \ar[u]
\\
\ar[u]
SH^*(M_+)^{<D^2\tau} \ar[r] & SH^*(M_-)^{<D^3\tau} \ar[ul] \ar[u]
\\
\ar[u]
SH^*(M_+)^{<\tau} \ar[r] & SH^*(M_-)^{<D\tau} \ar[ul] \ar[u]
}
\end{equation}

\begin{definition}
For a given Liouville domain $M$, let $r(M,\tau)$ be the total dimension of the image of the canonical map $SH^*(M)^{<\tau} \rightarrow SH^*(M)$. This is clearly an increasing function of $\tau$. Define the growth rate of $M$ to be
\begin{equation}
\Gamma(M) = \overline{\lim}_{\tau}\, \frac{r(M,\tau)}{\log(\tau)} \in \{-\infty\} \cup [0;\infty].
\end{equation}
\end{definition}

\begin{lemma}
$\Gamma(M)$ is invariant under Liouville isomorphism.
\end{lemma}

\proof Suppose that $SH^*(M_-) \iso SH^*(M_+)$ is infinite-dimensional, in which case both functions $r(M_{\pm},\tau)$ are unbounded by definition of direct limit. In particular, $\lim_{\tau \rightarrow \infty} \log(D)/r(M_+,\tau) = 0$. \eqref{eq:ladder} yields a chain of inequalities $r(M_+,\tau) \leq r(M_-,D\tau) \leq r(M_+,D^2\tau) \leq \cdots$, which imply that
\begin{equation}
\Gamma(M_+)^{-1} = \underline{\lim}_{\tau} \frac{\log(\tau)}{r(M_+,\tau)} =
\underline{\lim}_{\tau} \frac{\log(D\tau)}{r(M_+,\tau)}
\geq \underline{\lim}_{\tau} \frac{\log(D\tau)}{r(M_-,D\tau)} = \Gamma(M_-)^{-1}.
\end{equation}
The same holds in the other direction. In the remaining case, where the symplectic cohomology is finite-dimensional, the growth rates are either zero ($SH^* \neq 0$) or $-\infty$ ($SH^* = 0$). \qed

As before, the simplest situation is where the Reeb flow on $\partial M$ is a circle action. By a suitable choice of perturbation, one can then ensure that the number of generators in the chain complex $CF^*(H^\tau)$ grows linearly with $\tau$, which means that $\Gamma(M) \leq 1$. If no cancellation occurs, such as Examples \ref{th:surface} and \ref{th:cotangent-sphere}, equality holds; whereas for $M = D^{2n}$, we obviously still get $\Gamma(M) = -\infty$. In a different direction, one can look at cotangent bundles, where the issue is related to the more classical problem of growth rates of closed geodesics. For instance, the cotangent bundle of the torus has $\Gamma(D^*T^n) = n$, while the cotangent bundle of any hyperbolic manifold $N$ has $\Gamma(D^*N) = \infty$ due to a theorem of Margulis \cite{margulis69}.
Finally, we point out that a positive value of $\Gamma(M)$ gives a lower bound for the growth of the number of closed Reeb orbits on $\partial M$, for any generic choice of contact one-form (genericity is necessary because we want the orbits to be nondegenerate). For work in a somewhat similar spirit, see \cite{frauenfelder-schlenk06}.

\subsection{Affine varieties}
Let $X$ be a smooth projective variety, $\scrL \rightarrow X$ a holomorphic line bundle, and $s \in H^0(\scrL)$ a section vanishing along a divisor $D$. Choose a metric $||\cdot||$ on $\scrL$. Clearly,
\begin{equation} \label{eq:psi}
h = -\log ||s||
\end{equation}
is an exhausting function on the complement $U = X \setminus D$. Note that $dd^c h = 4i\,F|U$, where $F$ is the curvature of the unique connection on $\scrL$ compatible with the metric and holomorphic structure \cite[Chapter 1]{griffiths-harris}.

\begin{lemma} \label{th:finite-1}
Suppose that $D$ has normal crossings. Then the critical points of $h$ form a compact subset of $U$.
\end{lemma}

\proof If the statement is false, there must be a sequence of critical points converging to some limit in $D$.
Around the limit point, one can find local holomorphic coordinates $(x_1,\dots,x_n)$ and a local trivialization of $\scrL$, in which $s(x) = x_1^{w_1} \cdots x_k^{w_k}$ for some $1 \leq k \leq n$ and weights $w_1,\dots,w_k > 0$. Write the metric as $||\cdot|| = e^{\psi} |\cdot|$ with respect to that trivialization, so that
\begin{equation} \label{eq:d-h}
 dh = - d\psi - w_1\cdot d\log |x_1| - \cdots - w_k \cdot d\log |x_k|.
\end{equation}
If we take the vector field $Y = -x_1\partial_{x_1} -\cdots -x_k\partial_{x_k}$, then $Y.(\log |x_j|) = -1$ for $1 \leq j \leq k$, while $Y.\psi \rightarrow 0$ as $x \rightarrow 0$. Hence, $dh(Y)$ is bounded away from $0$ for small $|x|$, which is a contradiction. \qed

For the same reason:

\begin{lemma} \label{th:finite-2}
Assume again that $D$ has normal crossings. Suppose that we have two line bundles $\scrL_i$ ($i = 0,1$) with metrics $||\cdot||_i$ and sections $s_i$, such that $s_i^{-1}(0) = D$ (with possibly different multiplicities). Let $h_i$ be the associated functions. Then, the critical points of $(1-t)h_0 + th_1$, taken together for all $t \in [0;1]$, form a compact subset of $U$.
\end{lemma}

Hironaka's resolution theorem says that any smooth affine algebraic variety $U$ admits a good compactification. This means that $U \iso X \setminus D$, where $X$ is a smooth projective variety carrying an ample line bundle $\scrL$, and $D \subset X$ is a normal crossing divisor of the form $s^{-1}(0)$ for some holomorphic section $s \in H^0(\scrL)$. Ampleness means that there is a metric $||\cdot||$ on $\scrL$ for which $-iF$ is a K{\"a}hler form (in future, we will just say that the metric has positive curvature). This implies that $h$ is a plurisubharmonic function on $U$. Choose such a metric, and a $C$ which is larger than all the critical values of $h$ (this is possible by Lemma \ref{th:finite-1}). The corresponding sublevel set $M$ is a Liouville domain, just as in Example \ref{th:stein} (of which this is obviously a special case). The choice of $C$ is irrelevant up to deformation, hence up to Liouville isomorphism, see the discussion in Example \ref{th:unique-sublevel}. Similarly, Lemma \ref{th:finite-2} shows that $M$ is independent of the choice of metric, as well as the choice of line bundle $\scrL$. Independence of the compactification is slightly more tricky, but ultimately relies on the same kind of argument. Suppose that we have two good compactifications, of which one dominates the other. This means that there is a holomorphic map $\pi: X_0 \rightarrow X_1$, $D_0 = \pi^{-1}(D_1)$, which induces an isomorphism $X_0 \setminus D_0 \rightarrow X_1 \setminus D_1$. Compare $(\scrL_0,s_0)$ and the pullback $\pi^*(\scrL_1,s_1)$. The associated functions $h_0$, $\pi^*h_1$ are both plurisubharmonic on $X_0 \setminus D_0$, and Lemma \ref{th:finite-2} applies, providing a suitable family of Liouville structures interpolating between the two given ones. Another feature of resolution of singularities is that given two good compactifications, one can always find a third one which dominates the two.

We have just shown how to associate to any smooth affine variety a Liouville domain, in a way which is unique up to Liouville isomorphism. In particular, this allows us to speak of the symplectic cohomology of the affine variety itself, which we will denote by $SH^*(U)$, and of its growth $\Gamma(U)$.

\subsection{Local circle actions}
The issue of translating the insight provided by Hironaka's theorem into explicit control over the Reeb dynamics is elementary, but still tricky. We will only cover the simplest nontrivial case, that of algebraic surfaces.

\begin{theorem} \label{th:algebraic-surface}
For any smooth affine algebraic surface $U$, $\Gamma(U) \leq 2$.
\end{theorem}

As usual, we start by writing $U = X \setminus D$, where $D$ has normal crossings and is defined by a section $s$ of an ample line bundle $\scrL$. Since our definition of normal crossing divisor excludes self-intersections, the irreducible components $K \subset D$ are smooth.

{\bf Step 1:} Constructing compatible local circle actions. {\em For each $K$ we want to have a Hamiltonian circle action $\rho_K$, defined in a neighbourhood of $K$ inside $X$, which preserves $K$ and rotates its normal bundle. These actions should be mutually compatible, which means that near each crossing point $p \in K_1 \cap K_2$, $(\rho_{K_1},\rho_{K_2})$ form a Hamiltonian $T^2$-action.}

In the simplest case of a smooth $D = K$, this is quite elementary (it follows from the tubular neighbourhood theorem for symplectic submanifolds). However, compatibility for different components can only be achieved
if those components intersect orthogonally (with respect to the symplectic form). In the algebraic surface case, we can assume that $||\cdot||$ is chosen in such a way that around each singular point of $D$, we have local holomorphic coordinates in which $D = \{x_1x_2 = 0\}$, and where the K{\"a}hler form is standard. This is an easy result, proved by a local patching process (see \cite[Corollary 7.1]{ruan02b} or \cite[Lemma 1.7]{seidel01}; in fact, our Step 1 is generally quite similar to certain constructions in \cite{ruan02b}). Having made that choice, one takes both circle actions involved to be the standard linear ones in our local coordinates, and then extends them over the smooth parts of $D$ by a suitable choice of symplectic tubular neighbourhoods.

Throughout the rest of the argument, we will use a metric $||\cdot||$ of the kind constructed above, and the associated data (the K{\"a}hler form, the plurisubharmonic function $h$, the one-form $\theta = -d^ch$, and the Liouville vector field $Z = \nabla h$).

{\bf Step 2:} Making the Liouville flow symmetric near the crossings. {\em There is a smooth function $k$ on $X$ such that the one-form $\theta' = \theta + dk$ is $T^2$-invariant near each singular point of $D$.}

We work in the previously chosen local coordinates near such a singular point, and also choose a local trivialization of $\scrL$ in which $s(x_1,x_2) = x_1^{w_1}x_2^{w_2}$. By assumption, the metric $||\cdot|| = e^{\psi} |\cdot|$ gives rise to the standard K{\"a}hler form, and this implies that $d^c\psi + d^c(\quarter |x|^2)$ is a closed (hence exact) one-form. Take the function whose derivative is that one-form, and multiply it by a (radial) cutoff to make it zero outside a smal neighbourhood of the origin. This gives the desired $k$.

Let $Z'$ be the Liouville flow associated to $\theta'$. Near the origin, this is $-d^c(\quarter |x|^2 + |s(x)|)$, and by comparing that with \eqref{eq:d-h}, one sees that $Z'.h>0$. A similar argument applies near smooth points of $D$, using the fact that $Z'-Z$ is bounded (with respect to any metric on $X$). The consequence is that sublevel sets $M = \{h(x) \leq C\}$ ($C \gg 0$), equipped with $\theta'$, are Liouville domains. Moreover, Lemma \ref{th:moser} can be applied to show that this Liouville structure is isomorphic to the one obtained from the original $\theta$.

{\bf Step 3:} Making the boundary symmetric. {\em There are arbitrarily large subdomains $M' \subset U$, such that $Z'$ points strictly outwards along $\partial M'$, and such that $\partial M'$ is invariant under all the local circle actions $\rho_K$.}

Let's start by discussing the simplest case, which is when $K = D$ is smooth. Let $m$ be the moment map of the local circle action $\rho = \rho_K$, normalized so that $m|D = 0$. We claim that
\begin{equation} \label{eq:through}
 dm(Z') < 0
\end{equation}
at all points of $U$ which are sufficiently close to $D$. To see why that is true, fix a point of $D$, and choose local holomorphic coordinates around that point, in which $D = \{x_1 = 0\}$, and so that the K{\"a}hler form is standard at $x = 0$. This means that the induced circle action $D\rho$ on the tangent space at $x = 0$ is the standard rotation of the first component. Hence, $m(x) = \half |x_1|^2 + O(|x|^3)$, and $\nabla m = x_1 \partial_{x_1} + O(|x|^2)$. By arguing as in \eqref{eq:d-h} one sees that $dh(\nabla m) = dm(Z)$ is negative, and in fact bounded away from zero. The same holds if one replaces $Z$ with $Z'$, since the difference between the two is bounded. Having that, one would then define $M'$ to be the complement of the locus $\{m(x) < \epsilon\}$, for some small $\epsilon$ (recall that $m$ is defined only locally, which is why we formulate this in terms of complements).

Near the singularities of $D$ one uses a slightly more complicated local model. Consider the usual local coordinates, and let $m_1(x)= \half|x_1|^2$, $m_2(x) = \half|x_2|^2$ be the moment maps of the two commuting $S^1$-actions. Take a function $\kappa$ such that $\kappa(s_1,0) = \kappa(0,s_2) = 0$; and $\partial_{s_1}\kappa > 0$ if $s_2 > 0$, respectively $\partial_{s_2}\kappa > 0$ if $s_1 > 0$. One then takes $M'$ to be the complement of the subset where $\{\kappa(m_1(x),m_2(x)) < \epsilon\}$. Transversality of $Z'$ along $\partial M'$ is easy to prove in this case, since we have an explicit formula for the vector field close to $x = 0$.
By a suitable choice of $\kappa$, one can ensure that these local models can be patched into the previously described construction along the smooth part of $D$.

One can easily see from the construction that $Z'$ has only a compact set of stationary points. By taking $M'$ large enough, one can ensure that it contains all those points. In that case, a variation on the argument from Example \ref{th:unique-sublevel} shows that passing from $(M,\theta'|M)$ to $(M',\theta'|M')$ does not change the isomorphism type of the Liouville structure.

{\bf Step 4:} Making the Reeb flow symmetric. {\em One can find a contact one-form $\alpha'$ on $\partial M'$, whose exterior derivative is $\omega|\partial M' = d\theta'|\partial M'$, and which is invariant under the local $S^1$-actions.}

Again, we begin with the case where $K = D$ is smooth. In that situation, $Z'$ and $\rho_t^*Z'$, for any $t \in S^1$, both point outwards along $\partial M'$. Hence, one can average $\rho_t^*\theta'$, and obtains an $S^1$-invariant one-form defined near $\partial M'$, whose exterior derivative is $\omega$, and whose dual vector field points outwards. The restriction of that is the desired contact one-form. In the general case, where $D$ has normal crossings, the one-form $\theta'$ we have constructed is already invariant under the $T^2$-action near each crossing point. It is therefore sufficient to carry out the averaging process away from the crossings, which roughly speaking means along the smooth parts of $D$. One can apply a deformation (Gray's theorem) to prove that the new contact structure is isomorphic to the one given by $\theta'|\partial M'$. Hence, we can find an isomorphic Liouville structure on $M'$ whose boundary contact one-form is precisely this $\alpha'$.

The dynamics of the Reeb flow on $(\partial M',\alpha')$ can be explicitly described. Over the boundary parts lying close to the smooth points of $D$, the Reeb flow is a circle action. Near the singular points of $D$, the model is $\R \times T^2$, with a flow that translates $\{r\} \times T^2$ with some speed $\xi(r) \in (0;\infty)^2$. Here, $\xi(s) = (\xi_-,0)$ for $s \leq -S$, $\xi(s) = (0,\xi_+)$ for $s \geq S$ ($S \gg 0$), so that outside a compact subset one recovers the circle action. The precise form of $\xi$ depends on the choice of $\kappa$ in the construction, but one can ensure that on the interval $(-S;S)$, $\partial_s \xi_1 < 0$, $\partial_s \xi_2 > 0$. Then, an elementary argument shows that the number of periodic tori of our translational flow whose period is less than some given $\tau$ increases like $\tau^2$. After a standard Morse-Bott type perturbation argument, Theorem \ref{th:algebraic-surface} follows.

\begin{remark}
Another way of thinking of this description is that one starts with the standard contact structure on positive circle bundles over surfaces with boundary (the surfaces being the components of $D$ with neighbourhoods of the crossing points removed), and glues them together with intermediate $[-1;1] \times T^2$ pieces, on which the Reeb dynamics is the one described above. From this perspective, most of the work done above goes into showing that the outcome of this ``plumbing'' construction is actually contact isomorphic to the boundary of the Liouville domain defined using K{\"a}hler geometry.

There is also a different approach, pointed out to the author by Ivan Smith. Equip $X$ with a Lefschetz pencil whose fibre at infinity is $D$, and such that all other fibres meet the strata of $D$ transversally. This should give rise to an open book decomposition of $\partial M$, compatible with its contact structure in the sense of \cite{giroux02}. Moreover, the monodromy of this open book decomposition is completely reducible (has no pseudo-Anosov components). In particular, the number of essential periodic points of the monodromy (seen as a symplectic automorphism) grows at most linearly. One would then use the relation between monodromy and Reeb flow of a suitable contact one-form. This is more conceptual than our previous strategy, but there are still some potentially tricky details, such as the behaviour of the Reeb flow near the binding.
\end{remark}

Generally, one expects a bound $\Gamma(U) \leq n$ for all $n$-dimensional smooth affine varieties $U$ (a slightly sharper statement would be that the growth rate is bounded above by the complex codimension of the smallest stratum in $D$, where $U = X \setminus D$). The underlying intuitive picture of the Reeb flow is analogous to the one for surfaces, but the technical details are more involved. In a different direction, we should point out that there are other classes of Liouville domains which potentially can be analyzed by similar methods, such as Milnor fibres of isolated (smoothable) singularities; in that case, the geometry is governed ``from the inside'' (by the structure of the resolution of the singularity).

\section{Non-vanishing theorems\label{sec:nonzero}}

The symplectic cohomology of a Liouville domain also plays a role in the study of its (closed) Lagrangian submanifolds. One way of formulating this connection is to say that there is a natural map (of commutative rings; actually, even of BV algebras) from symplectic cohomology to the Hochschild cohomology of the Fukaya category \cite{seidel02}. In particular, if the Fukaya category is not empty, the symplectic cohomology must be nonzero. Note that a similar conclusion was already reached in \cite{viterbo97a}, in a somewhat different way (namely, by combining Theorem \ref{th:viterbo} and Viterbo functoriality). Here, we adopt a more elementary version of the approach in \cite{seidel02}, and focus on concrete applications, in particular the case of four-manifolds. This is (previously unpublished) joint work of Ivan Smith and the author.

\subsection{The basic construction\label{subsec:restrict}}
Let $M$ be a Liouville domain, and $L \subset M$ a closed Lagrangian submanifold which is exact, $[\theta|L] = 0 \in H^1(L;\R)$. To deal with the usual sign issues, we should also either assume that $L$ is $spin$, or otherwise take $\K = \Z/2$. In this situation, there is a canonical map
\begin{equation} \label{eq:quantum-restriction}
SH^*(M) \longrightarrow H^{*+n}(L;\K),
\end{equation}
which is constructed as follows. Take the definition of $SH^*(M) = HF^*(H)$ from \eqref{eq:sh-1}. Consider
solutions of Floer's equation \eqref{eq:floer} on a half-cylinder, with Lagrangian boundary conditions. This means
\begin{equation} \label{eq:semi-floer}
\left\{
\begin{aligned}
& u: [0;\infty) \times S^1 \longrightarrow \hat{M}, \\
& \partial_s u + J(\partial_t u - X) = 0, \\
& \textstyle\lim_{s \rightarrow +\infty} u(s,\cdot) = x(\cdot), \\
& u(0,\cdot) \in L.
\end{aligned}
\right.
\end{equation}
Due to the compactness and exactness of $L$, the action $A_H$ for any loop lying inside $L$ is bounded by some constant $C$, and one therefore gets an a priori bound for solutions of \eqref{eq:semi-floer}:
\begin{equation} \label{eq:semi-bound}
E(u) \leq -A_H(x) + C.
\end{equation}
By using the evaluation map at $(0,0)$ on the moduli space of solutions, one associates to each generator of the Floer complex $CF^*(H)$ a chain in $L$. \eqref{eq:quantum-restriction} is the induced map on cohomology (we have used the word ``chain'' loosely; there are several possible ways of making this rigorous, such as using pseudo-cycles, or else describing $H^*(L)$ by the Morse complex of an auxiliary function).

The composition of \eqref{eq:quantum-restriction} with \eqref{eq:classical} is the ordinary restriction map $H^*(M;\K) \rightarrow H^*(L;\K)$. The slickest way (in a formal TQFT sense) to see this is to introduce another equivalent definition of \eqref{eq:classical}. Namely, consider a continuation map equation \eqref{eq:cont} where at one end, $(H_-,J_-) = (H,J)$ are the data used in defining $SH^*(M)$, and at the other end $H_+ = 0$. Finite energy solutions of this equation extend smoothly over the puncture at $s = +\infty$, hence can be thought of as maps $\C \rightarrow \hat{M}$. Evaluation at that puncture then yields a map from locally finite chains on $\hat{M}$ (which are Poincar{\'e} dual to ordinary cochains) to Floer cochains, and \eqref{eq:classical} is the induced map on cohomology. As a consequence, its composition with \eqref{eq:quantum-restriction} is given by counting solutions of an inhomogeneous pseudo-holomorphic map equation on a closed disc $D$. One can deform the equation by turning off the inhomogeneous term, and then ends up with ordinary pseudo-holomorphic discs $(D,\partial D) \rightarrow (\hat{M},L)$. By the exactness assumption, all such discs are necessarily constant, which yields the desired result. Having that, one can draw the conclusion announced at the beginning of this section:

\begin{prop}[Viterbo] \label{th:no-exact}
If $M$ contains a closed exact Lagrangian submanifold $L$, then $SH^*(M) \neq 0$ (with coefficients $\K = \Z/2$, and also with arbitrary coefficients provided that $L$ is $spin$).
\end{prop}

As already observed in \cite{viterbo97a}, the exactness assumption can be weakened. For instance, it is enough to assume that (for some almost complex structure, which is of contact type at infinity) there are no non-constant holomorphic discs in $(\hat{M},L)$. Such generalizations of Proposition \ref{th:no-exact} are important in applications, and ideally, one would like the criteria on $L$ to be as flexible as possible (potentially, even allowing immersed Lagrangian submanifolds). While such a general theory still remains to be developed, one can get an idea of the issues involved by looking at a particularly simple case, namely that of surfaces in four-manifolds.

\subsection{Essential tori\label{subsec:essential}}
Let $M$ be a four-dimensional Liouville domain, and $L \subset M$ a Lagrangian torus. We drop the assumption that $L$ must be exact, but still require it to be Bohr-Sommerfeld, which means that $[\theta|L]$ is an integral class; this condition is just a technical simplification, and could be lifted if desired. We will work with coefficients in some field $\K \supset \Q$, which is of course no problem since $L$ is spin. Let $\Lambda = \K((t))$ be the field of Laurent series, $\Lambda^{\geq 0} = \K[[t]]$ the subring of actual power series, and $\Lambda^{>0} = t\K[[t]]$ its maximal ideal. Occasionally, rather than saying that some $z$ lies in $\Lambda^{>0}$, we will use the notation $z \in O(t)$.

Take an almost complex structure $J$ on $\hat{M}$, which as usual should be of contact type at infinity, and consider as before pseudo-holomorphic discs $u: (D,\partial D) \rightarrow (\hat{M},L)$. Each such disc has a Maslov index $\mu_L(u) \in 2\Z$, and the virtual dimension of the moduli space (of unparametrized discs) is $\mu_L(u) - 1$. In fact, due to the structure theorem of \cite{kwon-oh96} (see also \cite{lazzarini00}), a generic $J$ has the following properties:
\begin{equation} \label{eq:generic-j}
\parbox{32em}{\em There are no non-constant holomorphic discs of Maslov index $\leq 0$. All holomorphic discs of Maslov index $2$ are regular. Moreover, if we fix a class $\alpha \in H_1(L)$, the moduli space of Maslov index $2$ discs whose boundary represents that class is a compact oriented one-dimensional manifold.}
\end{equation}
Fix a $J$ which is generic in this sense; and denote by $k_\alpha(L,J) \in \Z$ the number of Maslov index $2$ discs going through a general point of $L$, with boundary class $\alpha \in H_1(L)$ (equivalently, the degree of the evaluation map, defined on the moduli spaces of discs with one marked boundary point). In addition, fix a class $a \in H^1(L;\Lambda^{>0})$, which we think of as a flat connection on the trivial $\Lambda$-bundle over $L$. The assumption $a \in O(t)$ ensures that the holonomy around a loop $\alpha$, namely $\exp(\int_\alpha a) \in 1 + \Lambda^{>0}$, makes sense. Define
\begin{equation}
\begin{aligned}
& m^0(L,J,a) = \sum_\alpha k_\alpha(L,J)\, t^{\int_\alpha \theta} \exp(\textstyle\int_\alpha a) \in \Lambda^{>0}, \\
& m^1(L,J,a) = \sum_\alpha k_\alpha(L,J)\, t^{\int_\alpha \theta} \exp(\textstyle\int_\alpha a)\, \alpha \in H_1(L;\Lambda^{>0}).
\end{aligned}
\end{equation}

We say that $L$ is Floer-theoretically essential if, for some choice of $J$ and $a$, $m^1(L,J,a) = 0$. The importance of this condition is that it allows one to extend part of the construction leading to \eqref{eq:quantum-restriction}. Pick a point $y \in L$, which is a regular value of the evaluation map on all the moduli spaces of $J$-holomorphic discs of Maslov index $2$. Assuming that $L$ is Floer-theoretically essential, we have an associated $a$ for which $m^1(L,J,a) = 0$. Consider solutions of the equation \eqref{eq:semi-floer} satisfying $u(0,0) = y$, and define a map $\epsilon: CF^*(H) \longrightarrow \Lambda$ by counting them with weights
\begin{equation} \label{eq:t-weight}
t^{\int_\alpha \theta} \exp(\textstyle\int_\alpha a).
\end{equation}
Here $\alpha(t) = u(0,-t)$ is the boundary loop in $L$ associated to $u$ (with its proper boundary orientation). To see the importance of the first term in \eqref{eq:t-weight}, note that instead of \eqref{eq:semi-bound} we now have an energy bound of the form $E(u) \leq \int_\alpha \theta - A_H(x) + C$, which ensures compactness for fixed $x$ and $\int_\alpha \theta$. The second term in \eqref{eq:t-weight} is needed to ensure that $\epsilon$ is a chain map. As usual, to prove that property, one analyzes the boundary points of one-dimensional moduli spaces. A priori, three kinds of points appear, which are drawn schematically in Figure \ref{fig:bubble}. Counting points of type (i) yields $\epsilon \circ \delta$; points of type (ii) actually do not appear, because the $J$-holomorphic bubble has virtual dimension $\geq 1$, leaving dimension $\leq -1$ for the main component; finally, the count points of type (iii) with a fixed main component always yields a multiple of $m^1(L,J,a) \cdot \beta$, for some class $\beta \in H_1(L)$ which is the boundary loop of the main component. Hence, the only nontrivial contribution comes from (i), ensuring that $\epsilon \circ \delta = 0$. By the same kind of argument, if $1 \in SH^*(M)$ is the image of the unit in ordinary cohomology under \eqref{eq:classical}, then $H(\epsilon)$ maps this element to the unit in $\Lambda$. Hence:
\begin{figure}[ht]
\begin{centering}
\includegraphics{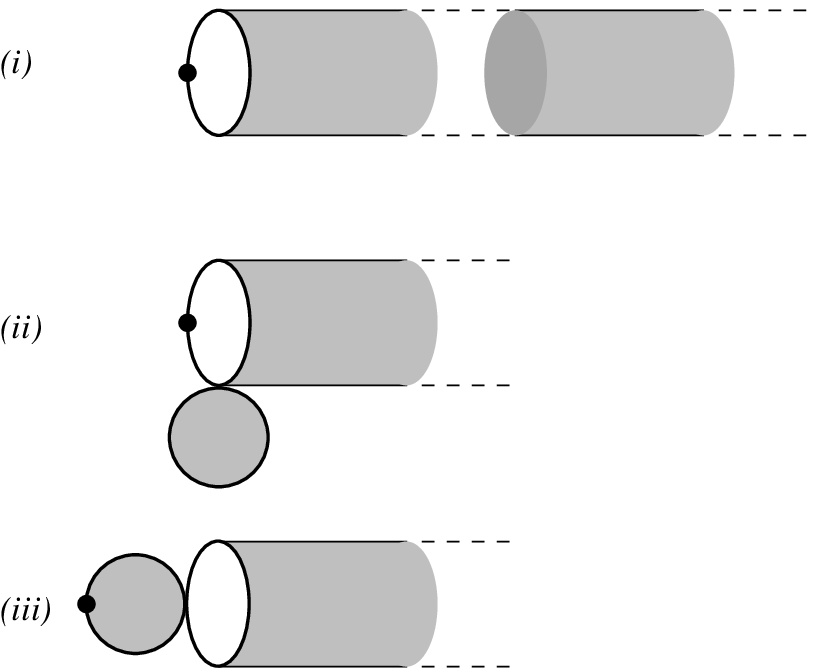}
\caption{\label{fig:bubble}}
\end{centering}
\end{figure}

\begin{prop}[Seidel-Smith] \label{th:no-essential}
If $M$ contains a Lagrangian torus which is (Bohr-Sommerfeld and) Floer-theoretically essential, then $SH^*(M) \neq 0$.
\end{prop}

\subsection{Wall-crossing}
So far, our definition of Floer-theoretically essential tori has skirted a crucial issue, which is dependence on the choice of almost complex structure. This is nontrivial because, inside the space of almost complex structures, there are codimension $1$ walls where Maslov index $0$ discs appear. Crossing such a wall will affect the numbers $k_\alpha(L,J)$, but one can compensate for that by a suitable modification of $a$. Unfortunately, because of the presence of multiple covers of Maslov index $0$ discs, it is hard to analyze the change in the moduli spaces explicitly. A more abstract approach, using virtual fundamental chains on suitable parametrized moduli spaces, was carried out in \cite{fooo} and yields the following result:

\begin{lemma}[Fukaya-Oh-Ohta-Ono] \label{th:fooo}
Suppose that, for some choice of $J$ which is generic in the sense of \eqref{eq:generic-j} and for some $c \in \Lambda^{>0}$, there exists an $a$ such that $m^0(L,J,a) = c$ and $m^1(L,J,a) = 0$. Then, the same is true for any other generic $J$ (with the same $c$, but generally a different $a$).
\end{lemma}

It is maybe helpful to review quickly the general formalism from \cite{fooo}. To every $L \subset M$, that formalism associates a curved (or obstructed) $A_\infty$-structure on $H^*(L;\Lambda^{\geq 0})$, which is given by a series of $t$-linear $\Z/2$-graded maps
\begin{equation} \label{eq:a-infty}
\mu^k : H^*(L;\Lambda^{\geq 0})^{\otimes k} \longrightarrow
H^{*+2-k}(L;\Lambda^{\geq 0}),
\end{equation}
$k \geq 0$, satisfying the generalized associativity equations. Specialization to $t = 0$ reproduces the classical $A_\infty$-structure on cohomology. The higher order terms are instanton corrections, coming from holomorphic discs. The $\mu^k$ are not unique, since they depend on many choices, both in the classical and instanton-deformed parts. However, given any two such structures $\{\mu^k\}$, $\{\tilde\mu^k\}$, one can find a curved $A_\infty$-homomorphism $\F$ between them, which for $t =0$ specializes to a quasi-isomorphism (whose linear part is the identity). Concretely, $\F$ is given by a sequence of maps
\begin{equation} \label{eq:functor}
\F^k: H^*(L;\Lambda^{\geq 0})^{\otimes k} \longrightarrow H^{*+1-k}(L;\Lambda^{\geq 0}),
\end{equation}
$k \geq 0$, satisfying the $A_\infty$-morphism equations, and with the additional property that $\F^0 \in O(t)$, $\F^1 \in id + O(t)$. Now consider odd degree elements $a \in H^*(L;\Lambda^{>0})$ satisfying the Maurer-Cartan equation
\begin{equation} \label{eq:mc}
\maurercartan(a) = \mu^0 + \mu^1(a) + \mu^2(a,a) + \cdots = 0.
\end{equation}
To each such element, one can associate a Floer cohomology group, which is the cohomology of the differential $d_a = \mu^1 + \mu^2(a,\cdot) + \mu^2(\cdot,a) + \cdots$ on $H^*(L;\Lambda)$. Using \eqref{eq:functor} one can transfer solutions of the Maurer-Cartan equation from one choice of $A_\infty$-structure to the other, and this induces isomorphisms of the associated Floer cohomology groups. For a more detailed exposition, see \S\S 7 (algebra) and 10 (geometry) of \cite{fooo}.

At this point, we want to specialize to $L = T^2$, and use the particular features of that situation to simplify the picture. First of all, the classical $A_\infty$-structure is known to be formal, so $\mu^k \in O(t)$ for $k \neq 2$, $\mu^2 = (\text{\it cup-product}) + O(t)$. Next, when computing the instanton corrections, one can choose an almost complex structure of type \eqref{eq:generic-j}. For our purposes, the most relevant part is to see what $\mu^k$ does to elements of $H^1$. For degree reasons, the only instanton correction comes from discs with Maslov index $2$, whose moduli spaces are smooth by assumption. A more careful computation shows that in fact,
\begin{equation}
\mu^k(a,\cdots,a) = \sum_{\alpha} k_{\alpha}(L,J) \, t^{\int_\alpha \theta} \, a(\alpha)^k/k! \in
\Lambda^{>0} = H^0(L;\Lambda^{>0}).
\end{equation}
Hence, solutions of \eqref{eq:mc} are precisely those $a$ such that $m^0(L,J,a) = 0$. Moreover, the associated Floer cohomology group is nontrivial if and only if $m^1(L,J,a) = 0$. The abstract theory ensures that the existence of such an $a$ is independent of all choices made in the construction.

This proves the $c = 0$ case of Lemma \ref{th:fooo}. Even though that case is actually sufficient for our applications, it is nevertheless instructive to briefly sketch the general argument, which is a little more sophisticated since it relies on the unitality on our $A_\infty$-structures. Specifically, one can arrange (by a suitable choice) that the unit $1 \in H^0(L;\Lambda^{\geq 0})$ satisfies $\mu^k(\cdots,1,\cdots) = 0$ for all $k \neq 2$, and $\mu^2(z,1) = z$, $\mu^2(1,z) = (-1)^{|z|} z$. Given any two choices of $A_\infty$-structures constructed in this particular way, the $A_\infty$-homomorphism from \eqref{eq:functor} will satisfy $\F^1(1) = 1$, and $\F^k(\cdots,1,\cdots) = 0$ for all $k \geq 2$. Instead of solutions of \eqref{eq:mc}, one can more generally consider projective Maurer-Cartan elements, which are those $a$ such that $\maurercartan(a) = c \cdot 1$. The resulting differentials $d_a$ still satisfy $d_a^2 = 0$, so there are well-defined Floer cohomology groups, and these are again independent of all choices involved.

\subsection{Crossing tori}
We now return to the situation discussed in Section \ref{sec:affine}, namely where $U = X \setminus D$ is an affine algebraic surface. Pick a singular point $p \in D$, and local holomorphic coordinates near that point, in which $D = \{x_1x_2 = 0\}$. As in the proof of Theorem \ref{th:algebraic-surface}, one can arrange that the K{\"a}hler form is standard in those coordinates. This means that for all sufficiently small $\epsilon_k > 0$,
\begin{equation} \label{eq:crossing-torus}
L = \{|x_1| = \epsilon_1, |x_2| = \epsilon_2\}
\end{equation}
is a Lagrangian torus in $U$, hence also in the associated Liouville domain $M$ (provided that this is chosen to be large enough). If the $\epsilon_k$ are rational multiples of $1/\pi$, the class $[\theta|L]$ is rational, and one can then rescale $\theta$ to make $L$ Bohr-Sommerfeld. We call such an $L$ a crossing torus associated to $p$.

Let $D' \subset D$ be a sub-divisor containing $p$, and denote by $D''$ the union of all components of $D$ which are disjoint from $D'$. Assume that the following conditions are satisfied:
\begin{itemize}
\item There is a neighbourhood $V \subset X$ of $D'$ containing $L$, such that the map $\pi_1(L) \rightarrow \pi_1(V \setminus D)$ is injective.
\item If $C$ is a holomorphic curve on $X$, whose intersection number with each component of $D''$ is zero or negative, then $C$ must be contained in $D'$.
\end{itemize}
If this holds, the crossing torus \eqref{eq:crossing-torus} is Floer-theoretically essential. The proof, which is a variation of an argument given in \cite[Section 5]{seidel-smith04b}, combines two techniques. First, using a elementary algebro-geometric trick, one degenerates $X$ to a singular variety (one component is the blowup of $X$ at $p$, and the other is $\CP{2}$). This degeneration actually yields a family of K{\"a}hler forms $\omega_t$ on $X$, together with a family of $\omega_t$-Lagrangian tori $L_t$. This family is parametrized by $t \in [0;\infty)$, with $t = 0$ being the original K{\"a}hler form and crossing torus, and moreover the class $[\omega_t] \in H^2(X,L_t)$ is constant throughout the deformation. By analyzing the limiting behaviour of pseudo-holomorphic discs, and using the geometric conditions above, one sees that if $t_k$ is a sequence going to $\infty$, and $u_k$ a sequence of non-constant holomorphic discs in $(X,L_{t_k})$, then the energy $\int_{u_k} \omega_k$ is necessarily unbounded. In informal language, there are no holomorphic discs in the limit $t \rightarrow \infty$. At that point, the proof of Lemma \ref{th:fooo} carries over with no changes, and establishes the desired property.

A concrete example is given by Ramanujam's surface \cite{ramanujam71}, which is an affine algebraic surface $U$ that is contractible (but not homeomorphic to $\R^4$, because its fundamental group at infinity is nontrivial). Proposition \ref{th:no-essential}, applied to the crossing torus introduced in \cite{seidel-smith04b}, shows that this surface has nonzero symplectic cohomology. Using the K{\"u}nneth formula from \cite{oancea04}, one sees that the same holds for $U \times U$ (which is simply-connected at infinity, hence actually diffeomorphic to $\R^8$, by the h-cobordism theorem). In particular, by looking at the associated Liouville domain, one gets the following consequence:

\begin{cor}[Seidel-Smith] \label{th:ramanujam}
There is a Liouville domain $M$ which is diffeomorphic to $D^8$ but has nonvanishing symplectic cohomology. In particular, $M$ is not Liouville isomorphic to standard $D^8$, and (in view of Corollary \ref{th:sphere} below) $\partial M$ is not contact isomorphic to standard $S^7$.
\end{cor}

\section{Handle attachment and classification problems\label{sec:novikov}}

We begin by reviewing another method of constructing Liouville domains, namely Weinstein handle attachment. This is inspired by the Morse theory of plurisubharmonic functions \cite{eliashberg90,eliashberg94}. Compared to Lefschetz fibrations, it is somewhat more flexible and economical, as well as closer to traditional manifold topology. As a cautionary note, we should point out that even this method has its limits: there are four-dimensional Liouville domains which have disconnected boundary, hence are not homotopy equivalent to $2$-dimensional cell complexes \cite{mcduff91b,geiges95}. Clearly, such a manifold cannot be decomposed into Weinstein handles (hence, is not a Weinstein manifold in the sense of \cite{eliashberg94}).

Attachment of $k$-handles with $k<n = \half \mathrm{dim} M$ is an essentially topological process, since it is governed by an h-principle \cite{eliashberg-mishachev}. On the side of pseudo-holomorphic curve theory, this is reflected by a remarkable result of Cieliebak, which says that the symplectic cohomology remains unchanged. In many cases, one can use $2$-handle attachment to kill the fundamental group, and this has implications for the classification of Liouville (and contact) structures.

\subsection{Weinstein handles\label{subsec:handles}}
Let $M$ be a Liouville domain, and $\xi = ker(\alpha)$ the contact hyperplane field on $\partial M$. Take a submanifold $\Sigma \subset M$ which is contact isotropic, so $\alpha|\Sigma = 0$. This means that that $T\Sigma$ is a subbundle of $\xi|\Sigma$, which is contained in its symplectic orthogonal complement $(T\Sigma)^{\perp,symp}$. The quotient $(T\Sigma)^{\perp,symp}/T\Sigma$ then inherits a symplectic structure; we call it the symplectic normal bundle of $\Sigma$, and denote it by $\nu^{symp}\Sigma$.

Now suppose that $\Sigma$ is a sphere of some dimension $k-1$, and that $\nu^{symp}\Sigma$ is trivial as a symplectic vector bundle. More precisely, we want to fix a diffeomorphism $f: S^{k-1} \rightarrow \Sigma$, or rather its class in $\pi_0(\mathrm{Diff}(S^{k-1},\Sigma)/O_k)$ (just as in Remark \ref{th:details}); and also a trivialization of $\nu^{symp}\Sigma$, at least up to homotopy. Then, one can form a new Liouville domain by attaching a Weinstein handle $H_k$ to the boundary of $M$ near $\Sigma$. Denote this operation by
\begin{equation} \label{eq:handle-attachment}
M' = M \cup_{\Sigma} H_k.
\end{equation}
$\partial M'$ is obtained from $\partial M$ by a surgery which replaces $\Sigma$ with a sphere of codimension $k$ (the new sphere is coisotropic, which means that its preimage in the symplectization $\R \times \partial M'$ is a coisotropic $\R \times S^{2n-k-1}$).

\begin{example} \label{th:1-handle}
Start with a disjoint union $M = M_1 \cup M_2$, fix one point $y_i$ in each $\partial M_i$, and then connect them by a one-handle. We will call the resulting $M'$ the boundary connected sum of the $M_i$, and denote it by $M_1 \#_{\partial} M_2$ (strictly speaking, this notation is appropriate only when $\partial M_i$ is connected). $\partial M'$ is the connected sum of $\partial M_1$ and $\partial M_2$ in the ordinary sense of the word. In this case, the choice of symplectic framing is irrelevant, since $\pi_0(Sp_{2n-2})$ is trivial; and indeed, the requirement that orientations should be preserved leaves no freedom in how to take the connected sum.
\end{example}

\begin{example} \label{th:c1-attach}
Let $M$ be a Liouville domain of dimension $\geq 6$, and $\Sigma \subset \partial M$ an isotropic loop which which represents the trivial class in $H_1(M)$. Attaching a handle creates a new generator $A \in H_2(M')$, which is unique up to adding classes in $H_2(M)$. The homotopy classes of symplectic framings form an affine space over $\pi_1(Sp_{2n-4}) \iso \Z$. Changing the framing by some integer raises the Chern number $\leftsc c_1(M'),A \rightsc$ by the corresponding amount.
\end{example}

As mentioned before, Weinstein handle attachment is derived from the Morse theory of plurisubharmonic functions. To build the local model, one starts with the quadratic function $h: \C^n \rightarrow \R$,
\begin{equation}
h(z_1,\dots,z_n) = \sum_{j=1}^{n-k} \quarter |z_j|^2 + \sum_{j=n-k+1}^n |re(z_j)|^2 - \half |im(z_j)|^2.
\end{equation}
The stable and unstable manifolds of the unique critical point $z = 0$ are the orthogonal subspaces $W^s = \{z_1 = \cdots = z_{n-k} = 0,\, z_{n-k+1} \in i\R, \dots, z_n \in i\R, h(z) = -1\}$, $W^u = \{z_{n-k+1} \in \R,\dots,z_n \in \R\}$. In relation to the previous picture, the original boundary $\partial M$ corresponds to the level set $\{h = -1\}$, and the sphere $\Sigma \subset \partial M$ to $W^s \cap \{h = -1\}$. Similarly, the new boundary would essentially be $\{h = 1\}$, containing the sphere $W^u \cap \{h = 1\}$.  We said ``essentially'' because the construction of the precise local model involves some cutoff functions, which change the level sets $\{h = -1\}$ and $\{h = 1\}$ to make them equal at infinity; see \cite{weinstein91}.

\begin{lemma} \label{th:cancel}
Suppose that $\Sigma$ bounds a contact isotropic $(k+1)$-disc $\Delta \subset \partial M$, in a way which is compatible with the symplectic framing of $\Sigma$. Then the handle attachment \eqref{eq:handle-attachment} can be cancelled. This means that there is a contact isotropic sphere $\Sigma' \subset \partial M'$ of dimension $k$, such that $M'' = M' \cup_{\Sigma'} H_{k+1}$ is Liouville isomorphic to $M$.
\end{lemma}

The meaning of the compatibility condition is as follows. Up to homotopy, there is a distinguished isomorphism $\nu^{symp}\Sigma \iso \nu^{symp}\Delta|\Sigma \oplus \C$. The unique trivialization of $\nu^{symp}\Delta$ induces one of $\nu^{symp}\Sigma$, and this must be homotopic to the given symplectic framing. We will not give a detailed proof of Lemma \ref{th:cancel}, but the idea can be explained in the basic model considered above. Inside $\{h = -1\}$ consider the contact isotropic submanifold given by
\begin{equation} \label{eq:boh}
z_1 \in \R^+,\, z_2 = \cdots = z_{n-k} = 0, \, z_{n-k+1} \in i\R, \dots, z_n \in i\R, h(z) = -1.
\end{equation}
This is diffeomorphic to $\R^+ \times S^{k-1}$, and its boundary is $W^s \cap \{h = -1\}$. In the situation of the Lemma, this corresponds to a neighbourhood of $\Sigma \subset \Delta$, hence closes up to a disc inside the whole of $\partial M$. On the other side $\{h = 1\}$, the same formula \eqref{eq:boh} defines a contact isotropic submanifold diffeomorphic to $\R^k$. Globally, this will close up to a contact isotropic sphere $\Sigma' \subset \partial M'$, which intersects the coisotropic $(2n-k-1)$-sphere coming from the previous surgery transversally in a single point. In this case, gluing a handle onto $\Sigma'$ will cancel the previous attachment (this is well-known in the differentiable category \cite{milnor65}, and is a folk theorem in the symplectic category).

In the situation considered in Example \ref{th:c1-attach}, any surface in $M$ bounding $\Sigma$ determines a preferred class $A \in H_2(M')$. If we take this surface to be a coisotropic disc $\Delta$ in $\partial M$, and attach the handle in a way which is compatible with the induced framing, $A$ will actually be represented by the isotropic two-sphere $\Sigma' \subset \partial M'$, which means that $\leftsc c_1(M'), A \rightsc = 0$.

\subsection{Subcritical handle attachment}
Let $M$ be a Liouville domain of dimension $2n$, and $\Sigma \subset \partial M$ a contact isotropic sphere of dimension $k-1<n-1$. Choose a symplectic framing, and perform the handle attachment \eqref{eq:handle-attachment}.

\begin{theorem}[Cieliebak] \label{th:cieliebak}
$SH^*(M') \iso SH^*(M)$.
\end{theorem}

This is Cieliebak's theorem in its original form \cite{cieliebak02}, which applies to $\Z/2$-graded symplectic cohomology groups. For our subsequent applications, we want $\Z$-gradings, hence have to sort out the relevant Chern class issues. If $k = 1$, any trivialization of the canonical bundle ${\mathcal K}$ can be extended from $M$ to $M'$. The extension is not necessarily unique, but Theorem \ref{th:cieliebak} holds equally for all choices. For $k = 2$, given a trivialization of ${\mathcal K}$ over $M$, any symplectic framing of $\Sigma$ determines a Maslov class $\mu_{\Sigma} \in H^1(\Sigma) \iso \Z$, which is the obstruction to extending the trivialization over $M'$. Changing the framings affects this class as in Lemma \ref{th:c1-attach}, so there is a unique choice which kills the obstruction. For higher $k$, there is no problem at all. We finish our discussion by mentioning a simple but beautiful application:

\begin{cor}[Smith] \label{th:sphere}
Let $M$ be a Liouville domain of dimension $2n \geq 4$, such that $\partial M$ is contact isomorphic to  standard $S^{2n-1}$. Then $SH^*(M) = 0$.
\end{cor}

\proof For $2n = 4$, Gromov \cite{gromov85} showed that $M$ must be Liouville isomorphic to standard $\R^4$, so there is nothing to prove. In higher dimension, we start with a theorem of Eliashberg-Floer-McDuff \cite{mcduff91b,eliashberg91}, which says that $M$ must be acyclic (and simply-connected, hence diffeomorphic to $D^{2n}$ by the h-cobordism theorem). In particular, $c_1(M) = 0$, so $\Z$-graded symplectic cohomology is well-defined. Using the spectral sequence \eqref{eq:morse-bott}, one immediately gets an upper bound $\mathrm{rank}\,SH^k(M) \leq 1$ for all $k \in \Z$. The same thing applies to the connected sum $M \#_{\partial} M$, but on the other hand we have $SH^*(M \#_{\partial} M) \iso SH^*(M) \oplus SH^*(M)$, which implies the result. \qed

\subsection{From groups to manifolds}
We will borrow liberally from Novikov's work on algorithmic decidability of problems in manifold topology. The original reference is \cite{novikov74}; for more recent expositions, see \cite[Chapter 2]{weinberger} (informal) or \cite{nabutovsky95} (more detailed). Let $P = \leftsc g_1,\dots,g_k\,|\,r_1,\dots,r_l \rightsc$ be a finite presentation of a group $\Gamma_P$, which satisfies $H_1(\Gamma_P) = H_2(\Gamma_P) = 0$. This class of groups is large enough so that the triviality problem, which is to decide given $P$ whether $\Gamma_P$ is trivial, is algorithmically unsolvable. Given any such $P$, one can algorithmically construct a homology $n$-sphere $S_P$ (for some fixed large $n$; $n = 6$ will do), such that $\pi_1(S_P) \iso \Gamma_P$. $S_P$ will come with distinguished loops $\gamma_j$ representing the generators $g_j$. Moreover, if the presentation is such that $\Gamma_P$ is the trivial group, then $S_P$ will actually be diffeomorphic to the standard sphere $S^n$.

We now symplectify this construction, killing the fundamental group in the process. After a small perturbation, we may assume that the loops $\gamma_j$ are embedded and disjoint. Choose a nowhere zero section of the (co)normal bundle of each $\gamma_j$. This defines an isotropic loop $\Sigma_j$ in the boundary of the cotangent bundle $D^*S_P$. Attach Weinstein two-handles to all the $\Sigma_j$, with the framing chosen in such a way that the resulting manifold $M_P$ has trivial first Chern class. It is easy to see that $M_P$ is simply-connected, and independent of all the choices made in the construction, up to Liouville isomorphism. Similarly, $\partial M_P$ is simply-connected, and independent of all choices up to contact isomorphism.
For purposes of comparison, we also choose a presentation $\bar{P}$ of the trivial group having the same number of generators as $P$, and form $M_{\bar{P}}$.

\begin{lemma}
$M_P$ is Liouville isomorphic to $M_{\bar{P}}$ if and only if $\Gamma_P$ is trivial.
\end{lemma}

\proof Theorems \ref{th:viterbo} and \ref{th:cieliebak} together show that
\begin{equation} \label{eq:viterbo-0}
SH^0(M_P) \iso SH^0(D^*S_P) \iso H_0({\mathcal L}S_P;\K)
\end{equation}
is generated by the conjugacy classes in $\Gamma_P$. Hence, its rank is $1$ iff $\Gamma_P$ is trivial. To prove the converse, suppose that $\Gamma_P$ is trivial. Then, both $M_P$ and $M_{\bar{P}}$ are constructed by starting with the same number of loops in $S^n$. One can deform one set of loops into the other (without crossings or self-intersections), and that induces an isotopy of the associated contact isotropic loops in $\partial D^*S^n$. \qed

\begin{lemma}
$\partial M_P$ is contact isomorphic to $\partial M_{\bar{P}}$ if and only if $\Gamma_P$ is trivial.
\end{lemma}

\proof One direction follows directly from the previous Lemma, but the other one needs a little more work.
Consider the trivial case $\bar{P}$, and write $\bar\gamma_j$ for the loops in $S_{\bar{P}} \iso S^n$ representing the generators.  Each such loop bounds a disc, and we can choose those discs to be embedded and mutually disjoint. One can then use the same method as before to produce contact isotropic discs $\bar\Delta_j$ whose boundary are the loops $\bar\Sigma_j$ used in the handle attachment process. From Remark \ref{th:c1-attach} and the discussion following Lemma \ref{th:cancel}, it follows that our choice of symplectic framings is compatible with $\bar\Delta_j$. Hence the handle attachment can be cancelled. In other words, there is a symplectic cobordism whose concave side is $\partial M_{\bar{P}}$, and whose convex side is $\partial D^*S^n$.

Now suppose that $\Gamma_P$ is nontrivial, but $\partial M_P \iso \partial M_{\bar{P}}$. In that case, we could attach $3$-handles to $M_P$ (equivalently, glue on the above-mentioned cobordism), so as to get another Liouville manifold $M$ whose boundary is (contact isomorphic to) $\partial D^*S^n$. By looking at Example \ref{th:cotangent-sphere}, one sees that the Reeb orbits on the boundary do not contribute any degree $0$ generators to symplectic cohomology, so $SH_0(M)$ must be one-dimensional. On the other hand, $SH_*(M) \iso SH_*(M_P)$ by Cieliebak's theorem, which contradicts \eqref{eq:viterbo-0}. \qed

\begin{cor} \label{th:unsolvable}
The following problems are algorithmically unsolvable: (i) the Liouville isomorphism problem for simply-connected Liouville domains; (ii) the computation of the rank of $SH^0(M)$ for a given simply-connected Liouville domain with $c_1(M) = 0$; (iii) the isomorphism problem for simply-connected closed contact manifolds.
\end{cor}

The discussion in Section \ref{sec:wild} will show that one can actually replace Liouville isomorphism in (i) with a more general notion (symplectic isomorphism of the completions).

In a different direction, one can restrict the class of presentations to ones with a fixed (sufficiently large) number of generators. The triviality problem remains unsolvable, and our argument shows that there is no argument which decides whether a given simply-connected contact manifold is isomorphic to $\partial M_{\bar{P}}$.
On the other hand, it is shown in \cite[Theorem 1]{nabutovsky-weinberger97} that the recognition problem is algorithmically solvable among the class of simply-connected smooth manifolds (of dimension $\geq 5$, up to diffeomorphism). Informally speaking, this means that in the world of classical high-dimensional manifold topology, the fundamental group is the only significant obstruction to detection problems. By comparing the two aspects, it follows that the problem of saying whether an arbitrary contact structure on the manifold $\partial M_{\bar{P}}$ is isomorphic to the given one is also algorithmically unsolvable. As an near-immediate consequence, that specific manifold must carry an infinity of non-isomorphic contact structures, which are moreover undistinguishable by classical invariants such as Chern classes. It is an interesting question to what extent one can simplify that ``undetectable'' manifold further, for instance, whether such a result also holds for the sphere.

Finally, on a more philosophical level, it's important to keep in mind the correct interpretation of results such as Corollary \ref{th:unsolvable}. The statement is that one cannot solve the isomorphism problem ``in a vacuum'', without having some specific insight or information about the manifolds concerned. One can argue that this situation rarely applies in practice (a possible exception are random constructions of the kind proposed in Section \ref{subsec:lefschetz}).

\section{Viterbo functoriality\label{sec:wild}}

Functoriality with respect to (a certain class of) embeddings is probably the most fundamental property of symplectic cohomology. The original motivation came from the case of cotangent bundles, where maps between free loop space homologies can be defined using generating function methods \cite{viterbo94}. This, and the subsequent Floer-theoretic reformulation \cite{viterbo97a}, have many immediate implications for embedding problems. Here, we pursue a far more modest aim, which is to use functoriality to make the framework of symplectic cohomology a little more general and flexible. We also include a simple but elegant application, discovered by Mark McLean \cite{mclean06}.

\subsection{Statement}
A Liouville embedding is an embedding $\iota$ of a Liouville domain $M$ into another such domain $M'$ (of equal dimension), such that $\iota^*\theta' = e^\rho\theta + d(\text{\it some function})$, for some $\rho \in \R$. Viterbo's construction \cite{viterbo97a} associates to each such embedding a pullback map
\begin{equation} \label{eq:pullback}
SH^*(\iota): SH^*(M') \longrightarrow SH^*(M).
\end{equation}
This is homotopy invariant (within the space of all such embeddings), and functorial (with respect to composition of embeddings). For instance, suppose that we enlarge $M$ to $M' = M \cup_{\partial M} ([0;R] \times \partial M)$, which means adding a finite piece of the completion. Then, the inclusion $M \rightarrow M'$ induces an isomorphism on $SH^*$, since it can be deformed to the identity by following the Liouville flow.
Another case when the pullback maps are isomorphisms is subcritical handle attachment. In fact, the proper statement of Theorem \ref{th:cieliebak} is that the map induced by inclusion yields an isomorphism $SH^*(M') \iso SH^*(M)$.

\subsection{Liouville manifolds}
Take a manifold $M$ (non-compact, and without boundary) which comes with a one-form $\theta$ such that $d\theta = \o$ is symplectic. We will call $M$ an exact symplectic manifold. An exact symplectic isomorphism between such manifolds is a diffeomorphism which preserves $\theta$ up to adding an exact one-form. We will say that $M$ is Liouville if there is a sequence of Liouville domains $M_k$ and embeddings $\iota_k: M_k \rightarrow M$, such that $\iota_k^*\theta = \theta_k + d(\text{\it some function})$, whose images exhaust $M$.
This condition is clearly invariant under exact symplectic isomorphisms. An obvious example is the completion of any Liouville domain. Another one is an arbitrary Stein manifold, with $\theta$ obtained from an exhausting plurisubharmonic function $h$ (unlike the situation in Example \ref{th:unique-sublevel}, we do not assume that the critical point set of $h$ is compact).

One defines the symplectic cohomology of a Liouville manifold to be
\begin{equation} \label{eq:infinite-sh}
SH^*(M) = \underleftarrow{lim}_k SH^*(M_k),
\end{equation}
where the connecting maps are induced by the inclusions via \eqref{eq:pullback}. Of course, this is not a correct definition in the abstract sense (we should take the inverse limit in the derived sense, which means on the level of chain complexes). However, for the purposes of the rather rudimentary discussion here, the difference does not matter. \eqref{eq:infinite-sh} is obviously independent of the choice of exhaustion, and is invariant under exact symplectic isomorphisms. Moreover, in the case of completions, it reproduces the symplectic cohomology of the original Liouville domain. In particular, it follows that the symplectic cohomology of a Liouville domain is in fact an (exact) symplectic invariant of its completion; which is a somewhat stronger statement than invariance under Liouville isomorphism.

\subsection{A finite type condition}
The following notion is due to Mark McLean: a Liouville manifold $M$ is called of finite type if, for some $k$ as in \eqref{eq:infinite-sh}, the canonical map $SH^*(M) \rightarrow SH^*(M_k)$ is injective. Of course, the same is then true for all $l \geq k$, and for similar reasons, this condition is independent of the choice of exhaustion. Clearly, any completion of a Liouville domain is of finite type. On the other hand, one can easily construct infinite genus surfaces which are of infinite type. Here is a more interesting case, taken from \cite{mclean06}, where the phenomenon is not topologically detectable:

\begin{theorem}[McLean]
There is a Liouville manifold $M$ diffeomorphic to $\R^8$, which is not of finite type.
\end{theorem}

To construct this, one starts with the Liouville domain from Corollary \ref{th:ramanujam}, here called $M_1$. In the next step, enlarge it slightly (by adding a finite piece of the cone), and then take the boundary connected sum of two copies of the result. The outcome is another Liouville domain $M_2$, which comes with a natural inclusion $M_1 \hookrightarrow M_2 \setminus \partial M_2$ (in fact, two such inclusions, but we only need one), compatible with the given one-forms. We repeat this construction inductively, and get an increasing sequence $M_1 \subset M_2 \subset M_3 \cdots$, whose union $M$ is a Liouville manifold. Topologically, each $M_k$ is diffeomorphic to $D^8$. Moreover,
$M_{k+1} \setminus (M_k \setminus \partial M_k)$ is a simply-connected h-cobordism from $\partial M_k \iso S^7$ to $\partial M_{k+1} \iso S^7$, hence diffeomorphic to $[0;1] \times S^7$. Therefore, $M$ itself is diffeomorphic to $\R^8$. On the other hand, we have
\begin{equation}
SH^*(M_{k+1}) \iso SH^*(M_k) \oplus SH^*(M_k),
\end{equation}
with the restriction map being projection to the first summand. Hence, the finite type condition is not satisfied.

As a consequence, we know that there are at least three symplectically distinct Liouville manifold structures on $\R^8$, namely the standard structure (symplectic cohomology vanishes), the completion of $M_1$ (symplectic cohomology is nonzero, but the manifold is of finite type), and $M$ (of infinite type). By comparison with Gompf's results about Stein structures on manifolds homeomorphic to $\R^4$ \cite{gompf98}, it seems quite possible that there is actually an uncountable family of infinite type Liouville structures on any $\R^{2n}$, $n>2$; so, this is another point where the established results fall far short of expectations.

\section{Algebraic structures\label{sec:operations}}

Symplectic cohomology carries a rich structure of operations. Unsurprisingly, these are defined by looking at Floer-type equations on more general Riemann surfaces. A particular feature of the symplectic cohomology situation is that there is an essential difference between inputs and outputs (positive and negative marked points on our surfaces), and that degenerations to nodal surfaces are allowed only under certain additional conditions. The closest relative of the resulting theory is string topology, but comparisons with classical (equivariant) cohomology and with symplectic field theory are also instructive.

\subsection{Riemann surfaces and operations}
We temporarily return to the model case of symplectizations $M = \R \times Y$, as in Section \ref{subsec:symplectizations}. Take $H$ of the form \eqref{eq:H-h}, and $J$ of contact type \eqref{eq:contact-type}. Let $\Sigma$ be a Riemann surface, carrying a real one-form $\beta$ which satisfies $d\beta \leq 0$ everywhere. Consider a solution of the equation
\begin{equation}
\label{eq:generalized-floer}
\left\{
\begin{aligned}
& u: \Sigma \longrightarrow M, \\
& (du-X \otimes \beta)^{0,1} = 0,
\end{aligned}
\right.
\end{equation}
and write $\rho = e^r \circ u$ as usual. Then (in any local coordinate $z = s+it$ on $\Sigma$)
\begin{equation}
\Delta \rho = |du - X \otimes \beta|^2 - \rho h''(\rho) \frac{d\rho \wedge \beta}{ds \wedge dt} -
 \rho h'(\rho) \frac{d\beta}{ds \wedge dt},
\end{equation}
so the maximum principle applies. Note that if we take $\Sigma = \R \times S^1$ and $\beta = dt$, \eqref{eq:generalized-floer} reduces to Floer's equation; and more generally, taking $\beta = f(s)\,dt$ yields a special case of the continuation map equation \eqref{eq:cont}.

Take the completion $\hat{M}$ of a Liouville domain. Equip it with an almost complex structure $J$ which is of contact type at infinity, and a Hamiltonian function $H$ of the kind used to define symplectic cohomology in \eqref{eq:sh-1}. Take a Riemann surface $\Sigma$ with $p > 0$ negative and $q \geq 0$ positive punctures, and choose tubular ends $(-\infty;0] \times S^1 \rightarrow \Sigma$, $[0;\infty) \times S^1 \rightarrow \Sigma$ near those punctures. Suppose that we have a $\beta \in \Omega^1(\Sigma)$ satisfying $d\beta \leq 0$ and which, on each strip-like end, is a positive constant multiple of $dt$. Then, by counting solutions of \eqref{eq:H-h} in $\hat{M}$ with appropriate asymptotics, one defines a map
\begin{equation} \label{eq:phi}
 \Phi_{\Sigma,p,q}: SH^*(M)^{\otimes q} \longrightarrow \big(SH^*(M)^{\otimes p}\big)[-n\chi(\Sigma)].
\end{equation}
For instance, $\Sigma = \C$ gives a canonical element $\Phi_{\C,1,0} \in SH^{-n}(M)$. This is in fact the image of the ordinary identity element $1 \in H^0(M;\K)$ under the map \eqref{eq:classical}; we have already encountered this construction in Section \ref{subsec:restrict}, where $\C$ was thought of as a partial compactification of $\R \times S^1$. Next, taking $\Sigma$ to be a three-punctured sphere (two positive, one negative) yields a product on $SH^*(M)$, which (by standard Floer-theoretic argument) is commutative and associative, and has $\Phi_{\C,1,0}$ as a unit. Note that as a consequence, $SH^*(M)$ vanishes if and only if \eqref{eq:classical} is zero.

As another instructive example, consider $\R \times S^1$ with both punctures negative. This means that one takes $\beta = f(s)dt$, where $f$ goes from positive ($s \ll 0$) to negative ($s \gg 0$). The resulting invariant is an element
\begin{equation} \label{eq:copairing}
\Phi_{\R \times S^1,2,0} \in SH^*(M)^{\otimes 2},
\end{equation}
which we can think of as a pairing on the dual space $SH^*(M)^\vee$. However, this pairing turns out to be highly degenerate. This is a consequence of a general property of our invariants \eqref{eq:phi}, namely that they allow degenerations of the Riemann surface to a nodal one, provided that the relevant one-forms can be chosen so that they vanish near the degeneration locus. In our case, this means that one can degenerate $\R \times S^1$ to a union of two discs intersecting at a point. As a consequence, \eqref{eq:copairing} is in fact the image of the diagonal class in $H^*(M;\K)^{\otimes 2} \iso H_{4n-*}(M \times M, \partial(M \times M);\K)$ under the map \eqref{eq:classical} (applied to both factors).

Finally, one can generalize the basic framework above by allowing families of Riemann surfaces with tubular ends. The result is a map like \eqref{eq:phi}, but with the degree shifted down by the dimension of the parameter space. The simplest example is to take $\R \times S^1$, with one positive and one negative puncture, and to take a one-parameter family in which one of the two tubular ends gets rotated once. The outcome is a map $\Delta: SH^*(M) \longrightarrow SH^{*-1}(M)$, called the BV (Batalin-Vilkovisky) operator. Next, the moduli space $\scrM$ of three-punctured spheres (two positive, one negative puncture) equipped with tubular ends is homotopy equivalent to $(S^1)^3$. Each class in $H_k(\scrM)$ gives rise to a map $SH^*(M)^{\otimes 2} \rightarrow SH^*(M)$ of degree $-n-k$. For $[point] \in H_0(\scrM)$, this is just the product mentioned above, and the remaining classes correspond to operations obtained by composing that product with $\Delta$. In fact, a classical result of Getzler \cite{getzler94} shows that the same is true for all moduli spaces of $(d+1)$-punctured spheres (one negative, $d$ positive punctures); the resulting operations are generated by $\Delta$ and the product, and equip $SH^*(M)$ with the structure of a BV algebra.

\subsection{$S^1$-equivariant symplectic cohomology} The chain complex underlying symplectic cohomology actually admits an infinite sequence of maps
\begin{equation}
\delta_k: CF^*(H) \longrightarrow CF^*(H)[1-2k], \;\; k \geq 0,
\end{equation}
of which the first one, $\delta_0 = \delta$, is the ordinary boundary operator. Recall that in practice, to define that operator one introduces a perturbation of $H$ which is $t$-dependent, hence breaks the $S^1$-symmetry in the original equation \eqref{eq:floer}. The higher order corrections $\delta_k$ are defined by looking at families of such perturbations, which concretely means that they count solutions to certain (parametrized) continuation map equations (see \cite{viterbo97b}, where this construction was introduced into Floer theory for the first time). The basic equations can be summarized by saying that the equivariant differential
\begin{equation}
\delta_{eq} = \delta_0 + u\delta_1 + u^2\delta_2 + \cdots,
\end{equation}
where $u$ is a formal variable of degree 2, squares to zero.  We take the $\K[[u]]$-module $K = \K((u))/u\K[[u]] = \bigoplus_{i \leq 0} \K u^i$, and define
\begin{equation}
SH^*_{eq}(M) = H^*(CF^*(H) \otimes K, \delta_{eq}).
\end{equation}
Note that with this definition, we get a bounded below exhausting $u$-adic filtration on our chain complex. The associated spectral sequence converges to $SH^*_{eq}(M)$; its $E_1$ term is
\begin{equation} \label{eq:u-ss}
E_1^{pq} = \begin{cases}
SH^{q-p}(M) & p \leq 0, \\
0 & p > 0,
\end{cases}
\end{equation}
and its first differential is the BV operator $\Delta$. In particular, if $SH^*(M)$ is acyclic, then so is the equivariant version. The same strategy can be used to prove equivariant analogues of other results concerning symplectic cohomology. Here, we only want to consider one particularly simple special case, namely when the Reeb flow on $\partial M$ is a free circle action. For simplicity, let's also take $\K \supset \Q$. In that case, we have an equivariant analogue of \eqref{eq:morse-bott}, namely a spectral sequence converging to $SH^*_{eq}(M)$, with starting term
\begin{equation} \label{eq:morse-bott-eq}
E_1^{pq} = \begin{cases}
H^{q+n}(M;K) & p = 0, \\
H^{p+q+n-1-p\mu}(\partial M/S^1;\K) & p< 0, \\
0 & p>0.
\end{cases}
\end{equation}

\begin{remark} \label{th:variants}
There are two other viable versions of the $S^1$-equivariant theory, namely: the one based on $CF^*(H)[[u]]$ (which is a completed tensor product, meaning that arbitrary formal series in the generators of $CF^*(H)$ are allowed); and its localized version, which allows finitely many powers of $u^{-1}$. We have chosen $K$-coefficients since the resulting complex is still relatively small (countably generated).

A word of caution is appropriate. In the classical framework of equivariant (de Rham, let's say) cohomology, the equivariant differential has only two terms, $\delta_{eq} = \delta_0 + u\delta_1$. This means that one can take coefficients in any $\K[u]$-module, which is more general than the setup above. This greater flexibility allows one to construct theories with better localization properties \cite{jones-petrack}. The price to pay is that the analogues of the spectral sequence in \eqref{eq:u-ss} generally do not converge, yielding more restrictive invariance properties (equivariant maps which are non-equivariant homotopy equivalences do not generally induce isomorphisms). It is not clear whether this setup can be applied meaningfully to symplectic geometry. Another instructive point of comparison is cyclic homology, see in particular the discussion in \cite{jones87}.
\end{remark}

\begin{example} \label{th:s2}
Take $M = D^*N$ with $N = S^2$. The equivariant version of Theorem \ref{th:viterbo} \cite{viterbo97b} says that in such a case, $SH^*_{eq}(M) = H_{-*}^{eq}(\scrL N)$, where equivariant homology of the free loop space is defined in the classical way (as the homology of the Borel construction). The spectral sequence \eqref{eq:morse-bott-eq} reproduces the Morse-Bott spectral sequence for the geodesic energy functional. We have
\begin{equation}
 E_1^{pq} = \begin{cases} \K & p = 0, q = 0, \\
 \K^2 & p = 0, q < 0 \text{ even}, \\
 \K & p<0, q = p \pm 1, \\
 0 & \text{otherwise.}
\end{cases}
\end{equation}
The $p = 0$ column contains a copy of $K \iso H_{-*}^{eq}(point;\K)$, and this clearly survives to $E_\infty$. By using that fact, as well as the action of $u$, and comparing the result with known computations of loop space homology \cite{hingston92}, one sees that precisely the differentials $d_r: E_r^{-r,-r-1} \rightarrow E_r^{0,-2r}$, $r \geq 1$, are nonzero. In this context we should also mention Goodwillie's theorem \cite{goodwillie85} which says that for any simply-connected $N$, the map $H^*_{eq}(\scrL N;\K) \rightarrow H^*_{eq}(point;\K)$ is an isomorphism modulo $u$-torsion.
\end{example}

The equivariant theory carries operations, which are of a somewhat different kind than before. To see the most basic instance of this, consider again the moduli space $\scrM$ of $3$-punctured discs with tubular ends. The associated operations are parametrized by equivariant homology classes of $\scrM$, with respect to the action of $(S^1)^3$ rotating the ends. More precisely, one should take equivariant homology with $\K[[u^*]]$-coefficients for the positive ends, and with $\K((u^*))/u^*\K[u^*]$-coefficients for the negative end ($u^*$ is dual to $u$, hence of degree $-2$; so the last-mentioned version is classical equivariant homology). This consists of a single copy of $\K$, located in degree $2$, and the resulting operation
\begin{equation} \label{eq:lie}
 SH^*_{eq}(M) \otimes SH^*_{eq}(M) \longrightarrow SH^*_{eq}(M)[n-2]
\end{equation}
turns out to be a graded Lie bracket on $SH^*_{eq}(M)[2-n]$ (more generally, one expects $SH^*_{eq}(M)$ to have the structure of a gravity algebra \cite{getzler94}, which is the ``equivariant analogue'' of a BV structure). There is also still a map $H^{*+n}(M;\K) \rightarrow SH^*_{eq}(M)$, but the image of this is now a degenerate ideal (the bracket of an element in the image with anything else vanishes).

\begin{remark}
At this point, it is tempting to speculate about the relation between symplectic cohomology and (relative) Gromov-Witten invariants (for related work concerning contact homology, see \cite{katz05}). As a relatively simple example, suppose that $X^{2n}$ is a closed symplectic manifold with integral symplectic class $[\omega]$, and with first Chern class $c_1 = 2[\omega]$. In that case, one expects to have a finite number $\Psi$ of ``lines'' (pseudo-holomorphic spheres with area $1$) going through a generic point. Now, suppose that $D \subset X$ is a symplectic hypersurface representing some multiple $k[\omega]$, and $M$ a suitable piece of the complement $X \setminus D$, which is a Liouville domain. \eqref{eq:morse-bott-eq} applies, with $\mu = 2-4/k$ (the fractional grading actually makes sense, but for our immediate considerations, we only need $k = 1,2$ anyway). In the case $k = 1$, part of the spectral sequence is a map
\begin{equation}
\begin{aligned}
& \K = H^0(D;\K) = H^0(\partial M/S^1;\K) = E_1^{-1,-n} \xrightarrow{d_1^{-1,-n}} \\ & \qquad \qquad \longrightarrow
E_1^{0,-n} = H^0(M;\K) \rightarrow H^0(M;\K) = \K.
\end{aligned}
\end{equation}
One can conjecture that this should be multiplication by $\Psi$. At least, this idea is compatible with the computations in Example \ref{th:s2}, where $\Psi = 2$, and where $d_1^{-1,-n}$ indeed turned out to be nonzero. Next, in the case $k = 2$, one expects $H^{*+n}(M;K)$ to survive to the $E_\infty$ term of the spectral sequence. Let's suppose for simplicity that the whole spectral sequence degenerates, and also ignore problems arising from the choice of splitting of the resulting filtration on $SH^*_{eq}(M)$. With that set aside, part of the Lie bracket \eqref{eq:lie} is a map
\begin{equation}
 \K = H^0(D;\K)^{\otimes 2} \subset SH^{1-n}_{eq}(M)^{\otimes 2} \longrightarrow
 SH^{-n}_{eq}(M) \rightarrow H^0(M;\K) = \K,
\end{equation}
which could potentially be multiplication by $2\Psi$. Unfortunately, for $k>2$ the situation becomes more complicated, since one apparently needs to pass to secondary (Massey product type) operations on $SH^*_{eq}(M)$. Still, it is possible that a suitably defined structure of chain level operations on symplectic cohomology does indeed contain $\Psi$, for all values of $k$.
\end{remark}


\end{document}